\documentclass[mathpazo]{cicp}

\usepackage{geometry}
\usepackage{amsmath}
\usepackage{amssymb}
\usepackage{amsthm}
\usepackage{bm}
\usepackage{graphicx}
\graphicspath{{./figures/}}
\usepackage{cite}
\usepackage{color}
\usepackage{tabu}
\usepackage{booktabs}
\usepackage{subfigure}
\usepackage{enumerate}
\usepackage{tabu}
\usepackage{multirow}
\usepackage{diagbox}
\newtheorem{rem}{Remark}

\newcommand*\diff{\mathop{}\!\mathrm{d}}

\begin{document}
\title{Asymptotic-Preserving Neural Networks for Multiscale Kinetic Equations}


\author[]{Shi Jin\affil{1,2, 4},
    Zheng Ma\affil{1,2,3,5}, and Keke Wu\affil{1}\comma\corrauth}
\address{\affilnum{1}\ School of Mathematical Sciences,
    Shanghai Jiao Tong University,
    Shanghai, 200240, P. R. China. \\
    \affilnum{2}\ Institute of Natural Sciences,
    MOE-LSC, Shanghai Jiao Tong University,
    Shanghai, 200240, P. R. China. \\
    \affilnum{3}\ Qing Yuan Research Institute,
    Shanghai Jiao Tong University,
    Shanghai, 200240, China.\\
    \affilnum{4}Shanghai Artificial Intelligence Laboratory, Shanghai, China.\\
    \affilnum{5}CMA-Shanghai, Shanghai Jiao Tong University, Shanghai, China}

\emails{{\tt wukekever@sjtu.edu.cn} (K.~Wu)}

\begin{abstract}
    In this paper, we present two novel Asymptotic-Preserving Neural Networks (APNNs) for tackling multiscale time-dependent kinetic problems, encompassing the linear transport equation and Bhatnagar-Gross-Krook (BGK) equation in all ranges of Knudsen number.  Our primary objective is to devise accurate APNN approaches for resolving multiscale kinetic equations, which is also efficient in the small Knudsen number regime.
    The first APNN for linear transport equation is based on even-odd decomposition, which relaxes the stringent conservation prerequisites while concurrently introducing an auxiliary deep neural network. We conclude that enforcing the initial condition for the linear transport equation with inflow boundary conditions is crucial for this network. For the Boltzmann-BGK equation, the APNN incorporates the conservation of mass, momentum, and total energy into the APNN framework as well as exact boundary conditions.
    A notable finding of this study is that approximating the zeroth, first, and second moments---which govern the conservation of density, momentum, and energy for the Boltzmann-BGK equation,  is simpler than the distribution itself. Another interesting phenomenon observed in the training process is that the convergence of density is swifter than that of momentum and energy.
    Finally, we investigate several benchmark problems to demonstrate the efficacy of our proposed APNN methods.
\end{abstract}


\maketitle

\section{Introduction}

In scientific modeling, kinetic equations describe the dynamics of density distribution of particles that collide between themselves, or interacting with a medium or external fields.
These equations are defined in the phase space, thus suffer from curse-of-dimensionality. In addition, they typically involve multiple spatial and/or temporal scales, characterized by the Knudsen number, as well as high dimensional nonlocal operators, hence present significant computational challenges in numerical simulations.
For a comprehensive overview, please refer to the literature sources~\cite{BGP,DP-Acta,jin2010asymptotic,weinan2011principles} for a review.

Deep learning methods and deep neural networks (DNNs) have garnered immense attention within the scientific community,
including the possibility of resolving partial differential equations (PDEs)~\cite{beck2020,cai2021least,E2018,liao2019deep,lyu2020mim,raissi2019physics,deepGalerkin2018,zang2020weak}.
To explore alternative machine learning approaches for solving partial differential equations, we refer {to} the exemplary review article~\cite{beck2020}.
The key motivation behind such methods is to parameterize the solutions or gradients of PDE problems using deep neural networks.
These methods ultimately culminate in a minimization problem that is typically high-dimensional and nonconvex. Unlike classical numerical methods, deep learning methods are mesh-free and can solve PDEs in high dimension, complex domains and geometries.
It is also advantageous to possess flexibility and ease of execution.
Nonetheless, deep learning methods have several potential drawbacks, including lengthy training times, a lack of convergence, and reduced accuracy.
The idea of operator learning, on the other hand, offers a method to resolve a class of PDEs by training the neural network just {\it once}~\cite{li2020fourier,lu2021learning,zhang2021mod,li2021physics,wang2021learning,xiong2023koopman}.
It is important to note, however, that a number of issues regarding the convergence theory remain unclear.

In recent years, there has been extensive research conducted on multiscale kinetic equations and hyperbolic systems by employing deep neural networks.
This research includes, but is not limited to the works cited in references~\cite{CLM,HJJL,wuAPNN,lu2022solving,li2022model,bertaglia2022asymptotic1,bertaglia2022asymptotic2,lou2021}.
In the design of DNNS, the definition of loss functions is crucial. There are numerous choices available to build the loss when given a PDE. For instance, the variational formulation (DRM), the least-squares formulation (PINN, DGM),
the weak formulation (WAN), etc.
Due to the presence of small scales, the vanilla Physics-Informed Neural Networks (PINNs) can perform poorly for resolving multiscale kinetic equations where small scales present~\cite{lu2022solving, wuAPNN}.
A natural question is what kind of loss is ``good''.
One important feature is to preserve important physical properties, such as conservation, symmetry, parity, entropy conditions, and asymptotic limits, etc.
In our previous work ~\cite{wuAPNN}, we developed a DNN for multiscale kinetic transport {equations} (with possible uncertainties) by creating a loss that can capture the limiting macroscopic behavior, as the Knudsen number approaches zero, satisfying a property known as Asymptotic-Preserving (AP), an asymptotic property known to be important in designing efficient numerical methods for multiscale kinetic equations \cite{jin2010asymptotic, jin2022Acta},
hence justifies the need to use Asymptotic-Preserving Neural Networks (APNNs).
This APNN method is based on the micro-macro decomposition, and we demonstrated that the loss is AP with respect to the Knudsen number when it tends to zero.
{It is worth noting that not all kinetic equations undergo a micro-macro decomposition. Even when applicable, designing neural networks that automatically satisfy conservation properties for nonlinear kinetic equations presents a challenging task \cite{lemou2008new,bennoune2008uniformly}. Taking the example of the nonlinear BGK equation, neural network approximations of solutions based on micro-macro decomposition need to adhere to the local conservation of mass, momentum, and energy for the collision operator. While this conservation can be incorporated as soft constraints in the loss function, as mentioned in our previous work, such an approach has proven to compromise precision\cite{wuAPNN}. Consequently, we introduced a mechanism for embedding a mass conservation constraint directly into the neural network\cite{wuAPNN}. However, achieving a neural network that simultaneously satisfies mass, momentum, and energy conservation for the nonlinear BGK equation remains a formidable challenge. }
    
In this work, we propose a different APNN that provides a framework for more general kinetic equations, with a generic way to handle both conservation and asymptotic properties in the loss functions. 
{We contemplate incorporating redundant constraints into the new AP loss function which deviates from conventional AP numerical schemes. 
This, however, may be not a requirement in classical AP numerical schemes.
The main idea is to involve both the original equation for macroscopic quantities, as well as the dynamic equations of the moments--that satisfy the conservation properties--in the loss functions.} 
Evolving both dynamics have been used earlier for multiscale kinetic equations that allow the preservations of important moments at ease \cite{filbet2010class, jin2011,  gamba2019}. 
When such an idea is built into the loss function of DNN, it helps to loosen the stringent prerequisites for conservation.
{For} the time-dependent linear transport equation, the proposed APNN technique relies on {an} even-odd decomposition.
The novel loss function exhibits uniform stability with respect to  the small Knudsen number, whereby the neural network solution formerly converges uniformly to the macro solution in the zero {Knudsen} number limit.
Furthermore, we propose an APNN for
the {\it nonlinear} Boltzmann-BGK equation, constructing deep neural networks that automatically satisfy the conservation of mass, momentum, and energy, akin to APNN based on a micro-macro decomposition \cite{bennoune2008uniformly} as in \cite{wuAPNN} poses a formidable challenge. 
{Building} the dynamics of moments into the loss function greatly {alleviates} this challenge.
Thus, we have advanced by developing an APNN methodology that unites the fundamental equation with the equation governing local conservation law for the nonlinear Boltzmann-BGK equation. 

An outline of this paper is as follows.
In Section 2, a {detailed} illustration of Asymptotic-Preserving Neural Networks for linear transport equation and Boltzmann-BGK equation and the construction of loss functions are given.
Numerous numerical examples are presented in Section 3 to demonstrate the effectiveness of the APNNs, as well as the issue of preserving positivity for density distribution in DNN.  While our numerical examples are one-dimensional in space, we do include multidimensional uncertain variables, in order to demonstrate the ability of the network for {high-dimensional} problems.
The paper is concluded in Section 4.

\section{Methodology}
There are three primary components to the DNN framework.
The initial component entails utilizing a neural network as an approximation to the solution.
The second component involves evaluating the difference between the approximate and exact solutions, which is achieved through population and empirical loss/risk.
Finally, the third component is an optimization algorithm that aids in locating a local minimum.

To solve partial differential equations using deep neural networks, the typical procedure is generally analogous:
\begin{enumerate}
    \item Modeling: define the loss/risk associated {with} a PDE;
    \item Architecture: build a deep neural network (function class) for the trail function;
    \item Optimization: minimize the loss over the parameter space.
\end{enumerate}

\begin{figure}[ht]
    \centering
    \includegraphics[width=0.45\textwidth]{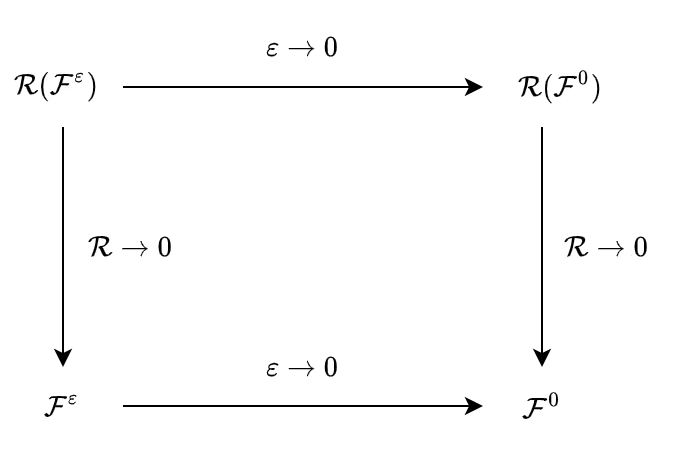}
    \caption{Illustration of APNNs.
        $\mathcal{F}^{\varepsilon}$ is the microscopic equation that depends on the small scale parameter $\varepsilon$ and $\mathcal{F}^{0}$ is its macroscopic limit as $\varepsilon \to 0$, which is independent of $\varepsilon$.
        The latent solution of $\mathcal{F}^{\varepsilon}$ is approximated by deep neural networks with its measure denoted by $\mathcal{R}(\mathcal{F}^{\varepsilon})$.
        The asymptotic limit of $\mathcal{R}(\mathcal{F}^{\varepsilon})$ as $\varepsilon \to 0$, if exists, is denoted by $\mathcal{R}(\mathcal{F}^{0})$. If $\mathcal{R}(\mathcal{F}^{0})$ is a good measure of $\mathcal{F}^{0}$, then it is called asymptotic-preserving (AP).}
    \label{fig:apnns}
\end{figure}

In terms of proposed APNNs, the primary element entails formulating a loss function that embodies the AP property~\cite{wuAPNN}.
The following diagram in Fig. \ref{fig:apnns} illustrates the idea of APNNs.

First, we introduce the conventional notations for DNNs\footnote{BAAI.2020.\ Suggested Notation
    for Machine Learning.\ https://github.com/mazhengcn/suggested-notation-for-machine-learning.}.
An $L$-layer feed forward neural network (or fully-connected neural network, FCNet) is defined recursively as,
\begin{equation}
    \begin{aligned}
        f_{\theta}^{[0]}(x) & = x,                                                                              \\
        f_{\theta}^{[l]}(x) & = \sigma \circ (W^{[l-1]} f_{\theta}^{[l-1]}(x) + b^{[l-1]}), \, 1 \le l \le L-1, \\
        f_{\theta}(x)       & = f_{\theta}^{[L]}(x) = W^{[L-1]} f_{\theta}^{[L-1]}(x) + b^{[L-1]},
    \end{aligned}
\end{equation}
where $W^{[l]} \in  \mathbb{R}^{m_{l+1}\times m_l}, b^{[l]}\in  \mathbb{R}^{m_{l+1}}, m_0, m_{L}$ are the input and output dimension, $\sigma$ is a scalar function and ``$\circ$'' means entry-wise operation.

Also, ResNet~\cite{he2016deep} is composed with several residual blocks with each part containing one input,
two weight layers, two activation functions, one identical (shortcut) connection,
and one output, which is defined recursively as,
\begin{equation}
    \begin{aligned}
        f_{\theta}^{[0]}(x) & = W^{[0]} x + b^{[0]},                                                                                                                                  \\
        f_{\theta}^{[l]}(x) & =  f_{\theta}^{[l-1]}(x) + \sigma \circ( W_2^{[l-1]} \sigma \circ (W_1^{[l-1]} f_{\theta}^{[l-1]}(x) + b_1^{[l-1]}) + b_2^{[l-1]}), \, 1 \le l \le L-1, \\
        f_{\theta}(x)       & = f_{\theta}^{[L]}(x) = W^{[L-1]} f_{\theta}^{[L-1]}(x) + b^{[L-1]},
    \end{aligned}
\end{equation}
here we use the same notations as previous for convenience.

We denote the set of parameters by $\theta$.
For simplicity of neural network presentation, we denote the layers by a list, i.e., $[m_0, \ldots, m_L]$.
For ResNet, the $i$-th block is denoted by $[m_i, m_i]$ and
one can find that $m_i (i = 1, \ldots, {L-1})$ is equal.

\subsection{APNN V2 based on even-odd decomposition}

Let us begin by examining the linear transport equation in the context of diffusive scaling, which can be expressed in the following form:
\begin{equation}\label{eqn: linear-transport-1d}
    \varepsilon \partial_t f + v\partial_x f = \frac{1}{\varepsilon}\left ( \frac{1}{2} \int_{-1}^1 f \, \diff v' - f \right ), \quad x_L < x < x_R, \quad -1 \leq v \leq 1,
\end{equation}
with in-flow boundary conditions as
{
\begin{equation}
F_{\text{B}} =
    \left \{
    \begin{aligned}
         & f(t, x_L, v) & = F_L(v), \quad \text{for} \quad v > 0, \\
         & f(t, x_R, v) & = F_R(v), \quad \text{for} \quad v < 0.
    \end{aligned}
    \right.
\end{equation}}
Here, $f(t, x, v)$ is the density distribution of particles at time $t \in \mathcal{T}$,
space point $x \in \mathcal{D}$, and traveling in direction $v \in \Omega := [-1, 1]$.
The parameter $\varepsilon > 0$ is the Knudsen number which denotes the ratio of the mean free path
over a characteristic length.
The initial function is given as a function of $x$ and $v$
\begin{equation}
    f(0, x, v) = f_0(x, v).
\end{equation}

By splitting equation and define even- and odd-parities as
\begin{equation}
    \begin{aligned}
        r(t, x, v) & = \frac{1}{2}[f(t, x, v) + f(t, x, -v)], \; 0 \le v \le 1,             \\
        j(t, x, v) & = \frac{1}{2\varepsilon}[f(t, x, v) - f(t, x, -v)],  \; 0 \le v \le 1,
    \end{aligned}
\end{equation}
one can obtain the following system of equations
\begin{equation}\label{eqn: even-odd}
    \begin{aligned}
         & \partial_t r + v\partial_x j = \frac{1}{\varepsilon^2}(\rho - r),                   \\
         & \partial_t j + \frac{1}{\varepsilon^2} v \partial_x r = -\frac{1}{\varepsilon^2} j,
    \end{aligned}
\end{equation}
where $\rho = \left \langle r \right \rangle := \int_0^1 r(t, x, v)  \diff v$.

Thus far, we have implemented the even-odd system, which is the building-block of an AP scheme for the problem \cite{JPT}. However,  simply employing neural networks for $r$ and $j$ does not result in an APNN framework. In order to establish an APNN framework, our next step involves the introduction of $\rho$ as a mediator between $r$ and $j$ within this system.

By integrating over $v$, the first equation gives
\begin{equation}\label{eqn: integral equation 1}
    \partial_t \left \langle r \right \rangle  + \int_0^1  v\partial_x j  \diff v = \frac{1}{\varepsilon^2} (\rho - \left \langle r \right \rangle),
\end{equation}
and since $\rho = \left \langle r \right \rangle$ one can write as follows
\begin{equation}\label{eqn: integral equation 2}
    \partial_t \rho + \int_0^1  v\partial_x j  \diff v = 0.
\end{equation}
Finally, Eq. (\ref{eqn: even-odd}) and Eq. (\ref{eqn: integral equation 2}) together with the constraint $\rho = \left \langle r \right \rangle$
constitute the even-odd formulation of Eq. (\ref{eqn: linear-transport-1d}) to be used for our APNN:
\begin{equation}\label{eqn: even-odd-system}
    \left \{
    \begin{aligned}
         & \varepsilon^2 \partial_t r + \varepsilon^2 v\partial_x j = \rho - r, \\
         & \varepsilon^2 \partial_t j + v \partial_x r = - j,                   \\
         & \partial_t \rho +  \left \langle v \partial_x j \right \rangle = 0,  \\
         & \rho = \left \langle r \right \rangle.
    \end{aligned}
    \right.
\end{equation}
When $\varepsilon \to 0$, the above equation formally approaches
\begin{equation}\label{eq:AP-limit}
    \left \{
    \begin{aligned}
         & r = \rho,                                                           \\
         & j = - v \partial_x r,                                               \\
         & \partial_t \rho +  \left \langle v \partial_x j \right \rangle = 0.
    \end{aligned}
    \right.
\end{equation}
Substituting the first equation into the second equation gives $j = - v \partial_x \rho$
and plugging into the third equation will result
\begin{equation}\label{eq:diffusion}
    \partial_t \rho -  \frac{1}{3} \partial_{xx} \rho = 0,
\end{equation}
which is exactly the diffusion equation.
\begin{rem}
    We emphasize that in Eq. (\ref{eqn: even-odd-system}),  the equation of local conservation law--or the moment equation
    $\partial_t \rho +  \left \langle v \partial_x j \right \rangle = 0$ is necessary in
    constructing the APNN loss in order to guarantee the AP property of the network. By coupling these equations of $r, j$ and $\rho$,
    one can obtain the loss for the diffusion limit equation when $\varepsilon \to 0$.
    {
    In ~\cite{wuAPNN}, we exemplify the AP numerical scheme based on the even-odd method to illustrate that not all classical AP numerical formats are suitable for constructing APNN methods. The even-odd formulation for the linear transport equation is  
    \begin{equation}\label{eq:even-odd}
        \begin{aligned}
             & \partial_t r + v \partial_x j                       = \frac{1}{\varepsilon^2}(\left \langle  r \right \rangle - r), \\
             & \partial_t j + \frac{1}{\varepsilon^2} v \partial_x r = -\frac{1}{\varepsilon^2} j.
        \end{aligned}
    \end{equation}
    It is noteworthy that within classical AP numerical schemes, the aforementioned Eq. (\ref{eq:even-odd}) is already sufficiently solvable.
    However if one uses it in the loss, it will not be an APNN. To see this, consider the case when $\varepsilon \to 0$. Since the DNN will only pick up the leading term, one has, in the $L^2$ sense,
    \begin{equation}
        \begin{aligned}
             & 0 = \left \langle  r \right \rangle - r,  \\
             & j = - v \partial_x r\,.
        \end{aligned}
    \end{equation}
    These two equations will not lead to the diffusion equation. This is because the corresponding loss function mentioned above does not provide the evolution process of the macroscopic quantity $\rho = \left \langle r \right \rangle$. This omission results in the neural network's approximated solution stagnating at the initial conditions.}
\end{rem}

For solving the linear transport equation by deep neural networks,
we need to use DNNs to parametrize three functions $\rho(t, x), r(t, x, v)$ and $j(t, x, v)$.
So here three networks are used. First,
\begin{equation}
    \rho^{\text{NN}}_{\theta}(t, x) := \exp \left( -\tilde{\rho}^{\text{NN}}_{\theta}(t, x)\right) \approx \rho(t, x),
\end{equation}

Second,
\begin{equation}
    r^{\text{NN}}_{\theta}(t, x, v) := \exp \left( -
    \frac{1}{2} (\tilde{r}^{\text{NN}}_{\theta}(t, x, v) + \tilde{r}^{\text{NN}}_{\theta}(t, x, -v) ) \right) \approx r(t, x, v),
\end{equation}
and
\begin{equation}
    j^{\text{NN}}_{\theta}(t, x, v) :=
    \tilde{\jmath}^{\text{NN}}_{\theta}(t, x, v) - \tilde{\jmath}^{\text{NN}}_{\theta}(t, x, -v)  \approx j(t, x, v),
\end{equation}
which automatically satisfy the even-odd properties.

Then we propose the least square of the residual of the even-odd system as the APNN loss

\begin{equation}\label{eq:loss-ap-even-odd}
    \mathcal{R}^{\varepsilon}_{\text{APNN}} = \mathcal{R}^{\varepsilon}_{\text{residual}} + \mathcal{R}^{\varepsilon}_{\text{constraint}} + \mathcal{R}^{\varepsilon}_{\text{initial}} + \mathcal{R}^{\varepsilon}_{\text{boundary}},
\end{equation}
{where $\mathcal{R}^{\varepsilon}_{\text{residual}}, \mathcal{R}^{\varepsilon}_{\text{constraint}}, \mathcal{R}^{\varepsilon}_{\text{initial}}, \mathcal{R}^{\varepsilon}_{\text{boundary}}$ are denoted by}
\begin{equation}\label{eq:loss-ap-even-odd-1}
    \begin{aligned}
        \mathcal{R}^{\varepsilon}_{\text{residual}} =   & \; \frac{\lambda_1}{|\mathcal{T} \times \mathcal{D} \times \Omega|} \int_{\mathcal{T}} \int_{\mathcal{D}} \int_{\Omega} | \varepsilon^2 \partial_t r^{\text{NN}}_{\theta} + \varepsilon^2 v\partial_x j^{\text{NN}}_{\theta} - (\rho^{\text{NN}}_{\theta} - r^{\text{NN}}_{\theta}) |^2 \diff{{v}} \diff{{x}}  \diff{t}                  \\
                                                        & + \frac{\lambda_2}{|\mathcal{T} \times \mathcal{D} \times \Omega|} \int_{\mathcal{T}} \int_{\mathcal{D}} \int_{\Omega} |\varepsilon^2 \partial_t j^{\text{NN}}_{\theta} + v \partial_x r^{\text{NN}}_{\theta} - (-j^{\text{NN}}_{\theta}) |^2 \diff{{v}} \diff{{x}}  \diff{t}                                                            \\
                                                        & + \frac{\lambda_3}{|\mathcal{T} \times \mathcal{D}|} \int_{\mathcal{T}} \int_{\mathcal{D}} | \partial_t \rho^{\text{NN}}_{\theta} +  \left \langle v\partial_x j^{\text{NN}}_{\theta} \right \rangle |^2   \diff{{x}}  \diff{t},                                                                                                         \\
        \mathcal{R}^{\varepsilon}_{\text{constraint}} = & \; \frac{\lambda_4}{|\mathcal{T} \times \mathcal{D}|} \int_{\mathcal{T}} \int_{\mathcal{D}} |\rho^{\text{NN}}_{\theta} -  \left \langle r^{\text{NN}}_{\theta} \right \rangle |^2   \diff{{x}}  \diff{t},                                                                                                                                \\
        \mathcal{R}^{\varepsilon}_{\text{initial}} =    & \; \frac{\lambda_5}{|\mathcal{D}|} \int_{\mathcal{D}} |\rho^{\text{NN}}_{\theta}(0, x) - \left \langle f_{0} \right \rangle |^2 \diff{\bm{x}} + \frac{\lambda_6}{|\mathcal{D} \times \Omega|} \int_{\mathcal{D}} \int_\Omega |\mathcal{I}(r^{\text{NN}}_{\theta} + \varepsilon j^{\text{NN}}_{\theta}) - f_{0}|^2 \diff{{v}} \diff{{x}}, \\
        \mathcal{R}^{\varepsilon}_{\text{boundary}} =   & \; \frac{\lambda_7}{|\mathcal{T} \times\partial \mathcal{D} \times \Omega|}  \int_{\mathcal{T}} \int_{\partial \mathcal{D}} \int_\Omega |\mathcal{B}(r^{\text{NN}}_{\theta} + \varepsilon j^{\text{NN}}_{\theta}) - F_{\text{B}}|^2 \diff{{v}} \diff{{x}} \diff{t}.
    \end{aligned}
\end{equation}
Here, $\lambda_i \; (i = 1, 2, \ldots, 7)$ are the penalty weights to be tuned and {$\mathcal{T}, \mathcal{D}, \Omega$ represent the bounded domains of time, space and velocity space, respectively.}
$|\mathcal{X}|$ denotes the measure of the domain $\mathcal{X}$.
Besides, $\mathcal{I}, \mathcal{B}$ are the initial and boundary operators.

Now the AP property of this loss can be carried out by considering its behavior for $\varepsilon$ small.
One may only need to focus on the first three terms of Eq. (\ref{eq:loss-ap-even-odd-1})
{
\begin{equation}
    \begin{aligned}
        \mathcal{R}^{\varepsilon}_{\text{residual}} = & \frac{\lambda_1}{|\mathcal{T} \times \mathcal{D} \times \Omega|} \int_{\mathcal{T}} \int_{\mathcal{D}} \int_{\Omega} | \varepsilon^2 \partial_t r^{\text{NN}}_{\theta} + \varepsilon^2 v\partial_x j^{\text{NN}}_{\theta} - (\rho^{\text{NN}}_{\theta} - r^{\text{NN}}_{\theta}) |^2 \diff{{v}} \diff{{x}}  \diff{t} \\
                                                            & + \frac{\lambda_2}{|\mathcal{T} \times \mathcal{D} \times \Omega|} \int_{\mathcal{T}} \int_{\mathcal{D}} \int_{\Omega} |\varepsilon^2 \partial_t j^{\text{NN}}_{\theta} + v \partial_x r^{\text{NN}}_{\theta} - (-j^{\text{NN}}_{\theta}) |^2 \diff{{v}} \diff{{x}}  \diff{t}                                        \\
                                                            & + \frac{\lambda_3}{|\mathcal{T} \times \mathcal{D}|} \int_{\mathcal{T}} \int_{\mathcal{D}} | \partial_t \rho^{\text{NN}}_{\theta} +  \left \langle v\partial_x j^{\text{NN}}_{\theta} \right \rangle |^2   \diff{{x}}  \diff{t}.
    \end{aligned}
\end{equation}
}
Sending $\varepsilon \to 0$, this will naturally lead to
{
\begin{equation}
    \begin{aligned}
        \mathcal{R}^{0}_{\text{residual}} = & \frac{\lambda_1}{|\mathcal{T} \times \mathcal{D} \times \Omega|} \int_{\mathcal{T}} \int_{\mathcal{D}} \int_{\Omega} |  \rho^{\text{NN}}_{\theta} - r^{\text{NN}}_{\theta} |^2 \diff{{v}} \diff{{x}}  \diff{t}                   \\
                                                  & + \frac{\lambda_2}{|\mathcal{T} \times \mathcal{D} \times \Omega|} \int_{\mathcal{T}} \int_{\mathcal{D}} \int_{\Omega} | v \partial_x r^{\text{NN}}_{\theta} - (-j^{\text{NN}}_{\theta}) |^2 \diff{{v}} \diff{{x}}  \diff{t}     \\
                                                  & + \frac{\lambda_3}{|\mathcal{T} \times \mathcal{D}|} \int_{\mathcal{T}} \int_{\mathcal{D}} | \partial_t \rho^{\text{NN}}_{\theta} +  \left \langle v\partial_x j^{\text{NN}}_{\theta} \right \rangle |^2   \diff{{x}}  \diff{t},
    \end{aligned}
\end{equation}
}
which is the least square loss of Eq.~\eqref{eq:AP-limit}
\begin{equation}
    \left \{
    \begin{aligned}
         & r = \rho,                                                           \\
         & j = - v \partial_x r,                                               \\
         & \partial_t \rho +  \left \langle v \partial_x j \right \rangle = 0.
    \end{aligned}
    \right.
\end{equation}
Same as previous derivation the third equation yields the diffusion Eq.~(\ref{eq:diffusion}).
Thus this proposed method is an APNN method.

\begin{rem}
{
    Additionally, one can try two Deep Neural Networks, namely $r_\theta$ and $j_\theta$, to approximate $r$ and $j$ without $\rho$ with the AP loss described by the following equations:
    \begin{equation}
        \begin{aligned}
             & \partial_t r + v \partial_x j = \frac{1}{\varepsilon^2}(\left \langle  r \right \rangle - r), \\
             & \partial_t j + \frac{1}{\varepsilon^2} v \partial_x r = -\frac{1}{\varepsilon^2} j, \\
             & \partial_t \left \langle r \right \rangle +  \left \langle v \partial_x j \right \rangle = 0.
        \end{aligned}
    \end{equation}
In our article, we use extra DNNs to approximate macroscopic quantities inspired by the relaxation schemes\cite{jin1995relaxation}.}
\end{rem}

\subsection{APNN for the Boltzmann-BGK equation}

The Boltzmann equation is the well-known kinetic model which captures the evolution of density distribution for rarefied gases~\cite{cercignani2000rarefied, kogan2013rarefied}.
The one-dimensional Boltzmann equation can be expressed in a dimensionless form as
\begin{equation}
    \partial_t f + v \cdot \nabla_x f = \frac{1}{\varepsilon} \mathcal{C}(f, f),  \quad t > 0, \quad (x, v) \in \mathbb{R} \times \mathbb{R},
\end{equation}
where the function, denoted by $f(t, x, v)$, is the density  distribution of particles,
while the right-hand side of the equation $\frac{1}{\varepsilon} \mathcal{C}(f, f)$ represents the term associated with bindary collisions between particles and is a non-linear operator.
Its action is generally limited to the velocity-dependent behavior of $f$ exclusively.
The parameter $\varepsilon > 0$ is the Knudsen number which denotes the ratio of the mean free path over a characteristic length.

The so-called Boltzmann collision operator $\mathcal{C}(f, f)$ possesses fundamental physical properties~\cite{bennoune2008uniformly} as
\begin{enumerate}
    \item conservation of mass, momentum and energy
          \begin{equation*}
              \int_{\mathbb{R}} m  \mathcal{C}(f, f) \mathrm{d} v = 0, \; m = {\left ( 1, v, \frac{1}{2} |v|^2\right )}^T.
          \end{equation*}
    \item The entropy dissipation inequality (the H-Theorem)
          \begin{equation*}
              \int_{\mathbb{R}} \mathcal{C}(f, f) \cdot \log (f) \mathrm{d} v \le 0.
          \end{equation*}
    \item The non-negative equilibrium functions f, namely those satisfying $\mathcal{C}(f, f) = 0$, correspond to the local Maxwellian distributions defined by
          \begin{equation}\label{eq:maxwellian}
              M(U) = \frac{\rho}{{(2 \pi T)}^{\frac{1}{2}}} \exp \left ( - \frac{|v - u|^2}{2 T}\right ),
          \end{equation}
          where density $\rho(t, x)$, macroscopic velocity $u(t, x)$ and temperature $T(t, x)$ of the gas are the continuum description by field variables defined by
          \begin{equation*}
              \rho = \int_\mathbb{R} f(v) \mathrm{d} v = \int_\mathbb{R} M(U) \mathrm{d} v,
          \end{equation*}
          \begin{equation*}
              u = \frac{1}{\rho} \int_\mathbb{R} v f(v) \mathrm{d} v = \frac{1}{\rho} \int_\mathbb{R} v M(U) \mathrm{d} v,
          \end{equation*}
          \begin{equation*}
              T = \frac{1}{\rho} \int_\mathbb{R} |u - v|^2 f(v) \mathrm{d} v = \frac{1}{\rho} \int_\mathbb{R} |u - v|^2 M(U) \mathrm{d} v.
          \end{equation*}

          Here, $U$ is the hydrodynamic variables (denisty, momentum and energy) which is a vector of the moments of $f$:
          \begin{equation}\label{eq:macro}
              U =
              \left (
              \rho,  \;
              \rho u, \;
              \frac{1}{2} \rho |u|^2 + \frac{1}{2} \rho T
              \right )^T = \int_\mathbb{R} m  f \mathrm{d} v,
          \end{equation}
          and density $\rho(t, x)$, macroscopic velocity $u(t, x)$ and temperature $T(t, x)$ are macroscopic quantities.
\end{enumerate}
As the frequency of collisions increases significantly, the mean free path, i.e., the distance a particle {travels} between two consecutive collisions,
becomes comparatively smaller than {the characteristic} length of the physical domain under consideration.
In such a scenario, a macroscopic portrayal of the gas seems more suitable.
The compressible Euler and compressible Navier-Stokes (CNS) equations serve as primary instances,
as they elucidate the dynamics of macroscopic quantities like the local density, momentum, and energy of the gas.
The CNS model surpasses the Euler equations in terms of accuracy, owing to its inclusion of factors such as viscosity and heat conductivity, which results in a correction of order $\varepsilon$.
Physically, classical fluid models may not fully capture the macroscopic evolution of gas, particularly when it is far away from {the} equilibrium state when the Knudsen number is large.

Fluid models, such as the compressible Euler or CNS type, are classically derived by using the moment method in conjunction with perturbation techniques like the Hilbert or Chapman-Enskog expansions~\cite{kogan2013rarefied}.
Specifically, the derivation of the CNS model from the Boltzmann equation in the fluid regime provides an approximation of viscosity and heat fluxes in the gas, up to the order of $\varepsilon^2$.
The primary challenge of solving the Boltzmann equation stems from the fact that the aforementioned term $\frac{1}{\varepsilon}$ stiffens as $\varepsilon$ approaches zero, thereby entering a fluid regime.
When considering this scenario, the resolution of the Boltzmann equation through a conventional explicit numerical method necessitates a time step of $\varepsilon$ magnitude.
This results in computationally costly operations when $\varepsilon$ is small.
At the level of Euler asymptotics, numerous authors have suggested asymptotically preserving numerical approximations to solve the Boltzmann equation.
For instance, numerical methods that are capable of capturing the accurate Euler limit have been proposed in~\cite{coron1991numerical, filbet2010class, bennoune2008uniformly}.

The intricate nature of the Boltzmann collision operator is circumvented by considering the simpler BGK model \cite{bhathnagor1954model} in our study as follows
\begin{equation}\label{eq:bgk}
    \partial_t f + v \cdot \nabla_x f = \frac{1}{\varepsilon} \left ( M(U) - f \right ),  \quad v \in \mathbb{R},
\end{equation}
here $f(t, x, v)$ is the density distribution of particles at time $t \in \mathcal{T}$,
space point $x \in \mathcal{D}$,
and traveling in direction $v \in \mathbb{R}$, $M(U)$ denotes the local Maxwellian distribution function.
Notice that the Boltzmann-BGK equation is an integro-differential equation with its nonlinear and non-local collision operator.
One can easily check that BGK operator satisfies mass, momentum and energy conservation.

Due to the properties of conserving mass, momentum and energy of collision operator,
one can multiply the BGK equation Eq. (\ref{eq:bgk}) by $m(v)$ and
then integrate them with respect to $v$ to obtain the following equations:
\begin{equation}\label{eq:claw1}
    \partial_t \left \langle m f \right \rangle + \nabla_x \cdot \left \langle v m f \right \rangle = 0, \; \text{where} \; \left \langle g \right \rangle =  \int_{\mathbb{R}} g(v) \mathrm{d} v,
\end{equation}
i.e.,
\begin{equation}\label{eq:claw2}
    \partial_t
    \begin{pmatrix}
        \rho   \\
        \rho u \\
        \frac{1}{2} \rho |u|^2 + \frac{1}{2} \rho T
    \end{pmatrix}
    + \nabla_x \cdot \left \langle v m f \right \rangle = 0.
\end{equation}
The equivalent formulation of Eq. (\ref{eq:claw2}) is
\begin{equation}\label{eq:fliud}
    \partial_t
    \begin{pmatrix}
        \rho   \\
        \rho u \\
        E
    \end{pmatrix}
    + \nabla_x \cdot
    \begin{pmatrix}
        \rho u       \\
        \rho u^2 + P \\
        E u + P u + Q
    \end{pmatrix}
    = 0,
\end{equation}
where $E = \frac{1}{2} \rho u^2 + \frac{1}{2} \rho T$ is the energy, $P = \int_\mathbb{R} |u - v|^2 f(v) \mathrm{d} v$ denotes the pressure tensor, and $Q = \int_\mathbb{R} (v - u) |u - v|^2 f(v) \mathrm{d} v$ represents the heat flux.
When taking the limit as $\varepsilon$ approaches zero, the density function $f$ will converge towards a local Maxwellian distribution $M(U)$. Under this approximation, the pressure tensor $P$ and the heat flux $Q$ can be expressed as $P = pI$ and $Q = 0$, respectively, where $p = \rho T$ represents the pressure and $I$ denotes the identity matrix.
Eq. (\ref{eq:fliud}) reduces to the compressible Euler equations
\begin{equation}\label{eq:fliud-limit}
    \partial_t
    \begin{pmatrix}
        \rho   \\
        \rho u \\
        E
    \end{pmatrix}
    + \nabla_x \cdot
    \begin{pmatrix}
        \rho u         \\
        \rho u^2 + p I \\
        (E + p) u
    \end{pmatrix}
    = 0.
\end{equation}

Finally, Eq. (\ref{eq:bgk}) and Eq. (\ref{eq:claw2}) with the constraint Eq. (\ref{eq:macro}) constitute the systems of BGK model for our APNN:
\begin{equation}\label{eqn:bgk}
    \left \{
    \begin{aligned}
         & \varepsilon \left ( \partial_t f + v \partial_x f \right ) = M(U) - f, \\
         & \partial_t U + \nabla_x \cdot \left \langle v m f \right \rangle = 0,  \\
         & U = \left \langle m f \right \rangle.
    \end{aligned}
    \right.
\end{equation}
The aforementioned idea was initially introduced in~\cite{jin2011}, wherein the employment of asymptotic-preserving schemes was demonstrated for the Fokker-Planck-Landau equation. Subsequently, it was further explored in~\cite{gamba2019}, focusing on the development of asymptotic-preserving numerical schemes and the conservation of numerical moments for collisional nonlinear kinetic equations.
While considering the Boltzmann-BGK equation, similar to APNN in \cite{wuAPNN},
it is crucial to highlight that our observation indicates the difficulty in creating a neural network for the non-equilibrium that adequately maintains the simultaneous conservation of mass, momentum and energy,
despite its micro-macro decomposition technique for resolution~\cite{bennoune2008uniformly}.
Thus, inspired by the idea of evolving both density distribution and moments for conservation properties~\cite{jin2011, gamba2019},
we incorporated the original equation concerning function $f$ into the system of local conservation laws Eq. (\ref{eq:claw1}) or Eq. (\ref{eq:claw2}),
thereby closing the newly established system of equations.

First we need to use DNNs to parametrize four functions $f(t, x, v), \rho(t, x), u(t, x)$ and $T(t, x)$.
So here four networks are used: $\rho_\theta(t, x), u_\theta(t, x), T_\theta(t, x), f_\theta(t, x, v)$. Besides, we restrict the range of velocity $\mathbb{R}$ to a bounded symmetrical domain $\Omega = [-V, V]$ and this assumption is used in most computations of the model.
The time and velocity variable $t, v$ are normalized into $[0, 1]$ and $[-1, 1]$ with scaling $\bar{t} = t / T, \bar{v} = v / V$.
It is worth pointing out that this normalization is necessary to alleviate the mismatching of the range of temporal domain and velocity domain.
For brevity we use the notation $W(t, x) = (\rho(t, x), \; u(t, x), \; T(t, x))^T$ and $W^{\text{NN}}_{\theta}(t, x)= (\rho^{\text{NN}}_{\theta}(t, x), \; u^{\text{NN}}_{\theta}(t, x), \; T^{\text{NN}}_{\theta}(t, x))^T$.

Then we propose the least square of the residual of BGK model as the APNN loss
\begin{equation}\label{eq:loss-ap-bgk}
    \mathcal{R}^{\varepsilon}_{\text{APNN}} = \mathcal{R}^{\varepsilon}_{\text{residual}} + \mathcal{R}^{\varepsilon}_{\text{claw}} + \mathcal{R}^{\varepsilon}_{\text{constraint}}
    + \mathcal{R}^{\varepsilon}_{\text{boundary}}
    + \mathcal{R}^{\varepsilon}_{\text{initial}},
\end{equation}
where $\mathcal{R}^{\varepsilon}_{\text{residual}}, \mathcal{R}^{\varepsilon}_{\text{claw}}, \mathcal{R}^{\varepsilon}_{\text{constraint}},
    \mathcal{R}^{\varepsilon}_{\text{boundary}},
    \mathcal{R}^{\varepsilon}_{\text{initial}}$ are denoted by
{
\begin{small}
    \begin{equation}\label{eq:loss-ap-bgk-part}
        \begin{aligned}
            \mathcal{R}^{\varepsilon}_{\text{residual}}   & =  \frac{\lambda_1}{| \mathcal{T} \times \mathcal{D} \times \Omega |} \int_{\mathcal{T}} \int_{\mathcal{D}} \int_{\Omega}| \varepsilon (\partial_t f^{\text{NN}}_{\theta} +  v\nabla_x f^{\text{NN}}_{\theta}) - \left ( M(U^{\text{NN}}_{\theta}) - f^{\text{NN}}_{\theta} \right ) |^2 \diff{{v}} \diff{{x}}  \diff{t}, \\
            \mathcal{R}^{\varepsilon}_{\text{claw}}       & = \frac{{\lambda_2}}{|\mathcal{T} \times \mathcal{D}|} \int_{\mathcal{T}} \int_{\mathcal{D}} |\partial_t U^{\text{NN}}_{\theta} + \nabla_x \left \langle v m f^{\text{NN}}_{\theta} \right \rangle |^2 \diff{{x}}  \diff{t},                                                                                        \\
            \mathcal{R}^{\varepsilon}_{\text{constraint}} & =  \frac{{\lambda_3}}{|\mathcal{T} \times \mathcal{D}|} \int_{\mathcal{T}} \int_{\mathcal{D}}  |U^{\text{NN}}_{\theta} - \left \langle m f^{\text{NN}}_{\theta} \right \rangle |^2 \diff{{x}}  \diff{t},                                                                                                            \\
            \mathcal{R}^{\varepsilon}_{\text{boundary}}   & = \frac{{\lambda_4}}{|\mathcal{T} \times \partial D|} \int_{\mathcal{T} \times \partial D} |W^{\text{NN}}_{\theta}(t, x) - W(t, x)|^2 \diff{x} \diff{t},                                                                                                                                                                  \\
            \mathcal{R}^{\varepsilon}_{\text{initial}}    & = \frac{{\lambda_5}}{\mathcal{|D|}} \int_{\mathcal{D}} |W^{\text{NN}}_{\theta}(0, x) - W(0, x)|^2 \diff{x} + \frac{\lambda_{6}}{|\mathcal{D} \times \Omega |} \int_{\mathcal{D}} \int_{\Omega}  |f^{\text{NN}}_{\theta}(0, x, v) - f(0, x, v)|^2 \diff{v} \diff{x},
        \end{aligned}
    \end{equation}
\end{small}
}
and $U^{\text{NN}}_{\theta}, W^{\text{NN}}_{\theta}$ are computed by $\rho^{\text{NN}}_{\theta}, u^{\text{NN}}_{\theta}$ and $T^{\text{NN}}_{\theta}$.

The AP property of this loss can now be realized by examining its behavior for infinitesimally small values of $\varepsilon$. One may only need to focus on the first two terms of Eq. (\ref{eq:loss-ap-bgk})
\begin{equation}
    \mathcal{R}^{\varepsilon}_{\text{APNN}} = \mathcal{R}^{\varepsilon}_{\text{residual}} + \mathcal{R}^{\varepsilon}_{\text{claw}},
\end{equation}
{
\begin{small}
    \begin{equation}
        \begin{aligned}
            \mathcal{R}^{\varepsilon}_{\text{residual}} & =  \frac{\lambda_1}{| \mathcal{T} \times \mathcal{D} \times \Omega |} \int_{\mathcal{T}} \int_{\mathcal{D}} \int_{\Omega}| \varepsilon (\partial_t f^{\text{NN}}_{\theta} +  v\nabla_x f^{\text{NN}}_{\theta}) - \left ( M(U^{\text{NN}}_{\theta}) - f^{\text{NN}}_{\theta} \right ) |^2 \diff{{v}} \diff{{x}}  \diff{t}, \\
            \mathcal{R}^{\varepsilon}_{\text{claw}}     & = \frac{{\lambda_2}}{|\mathcal{T} \times \mathcal{D}|} \int_{\mathcal{T}} \int_{\mathcal{D}} |\partial_t U^{\text{NN}}_{\theta} + \nabla_x \left \langle v m f^{\text{NN}}_{\theta} \right \rangle |^2 \diff{{x}}  \diff{t}.
        \end{aligned}
    \end{equation}
\end{small}
}
Sending $\varepsilon \to 0$, this will naturally lead to
{
\begin{small}
    \begin{equation}
        \begin{aligned}
            \mathcal{R}^{0}_{\text{residual}} & =  \frac{\lambda_1}{| \mathcal{T} \times \mathcal{D} \times \Omega |} \int_{\mathcal{T}} \int_{\mathcal{D}} \int_{\Omega}| M(U^{\text{NN}}_{\theta}) - f^{\text{NN}}_{\theta} |^2 \diff{{v}} \diff{{x}}  \diff{t},                 \\
            \mathcal{R}^{0}_{\text{claw}}     & = \frac{{\lambda_2}}{|\mathcal{T} \times \mathcal{D}|} \int_{\mathcal{T}} \int_{\mathcal{D}} |\partial_t U^{\text{NN}}_{\theta} + \nabla_x \left \langle v m f^{\text{NN}}_{\theta} \right \rangle |^2 \diff{{x}}  \diff{t},
        \end{aligned}
    \end{equation}
\end{small}
}
which is the least square loss of fluid dynamics equations, i.e., Eq. (\ref{eq:fliud-limit}).

\begin{figure}[ht]
    \centering
    \includegraphics[width=\textwidth]{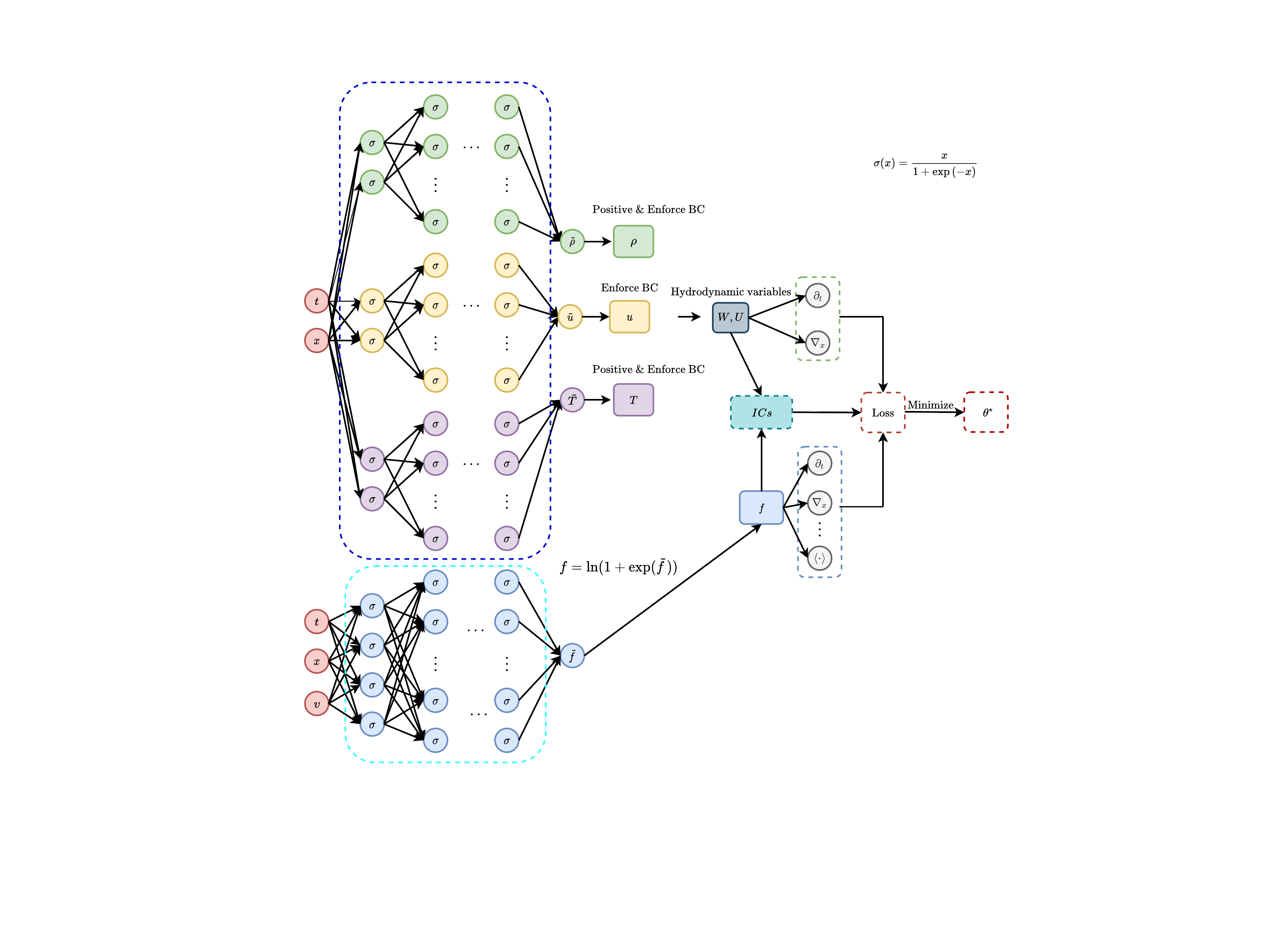}
    \caption{Schematic of APNNs for solving the Boltzmann-BGK equation.}\label{fig:APNNs}
\end{figure}

Thus, we have check that this loss is AP by sending $\varepsilon \to 0$.
Finally we put a schematic plot of our APNN method for Boltzmann-BGK problem in Fig. \ref{fig:APNNs}.

\section{Numerical results}

In this section, some numerical experiments will be  carried out to verify the performance of our proposed APNN methods.
Since the operator of $v$ in the loss of APNNs are integrals,
we approximate them with the Gauss-Legendre quadrature rule.

The reference solutions are obtained by standard finite difference method and
we will check the relative $\ell^2$ error of the solution $s(x)(\rho, u, T)$ of the APNN method, e.g.\ for $1$d case,

\begin{equation}
    \text{error} := \sqrt{
    \frac{\sum_j |s^{\text{NN}}_{\theta, j} - s^{\text{ref}}_j|^2}
    {\sum_j |s^{\text{ref}}_j|^2}
    }.
\end{equation}

\subsection{Experiment setting}

The activation function we used is $\sigma(x) = \tanh(x)$ for linear transport problems
and $\sigma(x) = x / (1 + \exp(-x))$ for BGK problems.
More specifically, fully-connected neural network is applied {to} the linear transport problem
while ResNet is used for the BGK problem for better performance.
In particular, the spatial domain of interest, denoted as $\mathcal{D}$, {includes} the interval $[0, 1]$ for the linear transport equation and $[-0.5, 0.5]$ for the Boltzmann-BGK equation. Additionally, the parameter $V$ is assigned a value of 10.
The number of quadrature points is $30$ for the linear transport equation and $64$ for {the} Boltzmann-BGK equation.
To train the networks, the Adam~\cite{kingma2014adam} version of the gradient descent methods is used to solve the optimization problem with Xavier initialization.
In practice of these cases, we need to tune the hyperparameters, such as neural network architecture, learning rate, batch size and so on,
to obtain a good level of accuracy~\cite{XDE}. As a matter of experience one may tune the weights of loss terms to make them at the same level
and a decreasing annealing schedule for learning rate is used for better numerical performance.
We use an exponential decay strategy for an initial learning rate $\eta_0 = 10^{-3}$ with a decay rate of $\gamma = 0.96$ and a decay step of $p = 200$ iterations:
{
\begin{equation*}
    \eta_i = \eta_0 \cdot \gamma^{\lfloor \frac{i}{p} \rfloor},
\end{equation*}
here, the variable $i$ represents the current $i-$th iteration step, and the symbol $\lfloor \cdot \rfloor$ denotes the floor function.}

\subsection{Problem 1: APNNs for solving linear transport equations}

Consider the linear transport equation:
\begin{equation*}
    \varepsilon \partial_t f + v\partial_x f = \frac{1}{\varepsilon}\left ( \frac{1}{2} \int_{-1}^1 f \, \diff v' - f \right ),
\end{equation*}
and recall that the even-odd system of the equation is
\begin{equation*}
    \left \{
    \begin{aligned}
         & \varepsilon^2 \partial_t r + \varepsilon^2 v\partial_x j = \rho - r, \\
         & \varepsilon^2 \partial_t j + v \partial_x r = - j,                   \\
         & \partial_t \rho +  \left \langle v \partial_x j \right \rangle = 0,  \\
         & \rho = \left \langle r \right \rangle.                               \\
    \end{aligned}
    \right.
\end{equation*}

The APNN empirical risk for the even-odd system of linear transport equation is
\begin{equation}
    \mathcal{R}^{\varepsilon}_{\text{APNN, transport}} = \mathcal{R}^{\varepsilon}_{\text{residual}} + \mathcal{R}^{\varepsilon}_{\text{constraint}} + \mathcal{R}^{\varepsilon}_{\text{initial}} + \mathcal{R}^{\varepsilon}_{\text{boundary}},
\end{equation}
{where $\mathcal{R}^{\varepsilon}_{\text{residual}}, \mathcal{R}^{\varepsilon}_{\text{constraint}},
\mathcal{R}^{\varepsilon}_{\text{initial}},
\mathcal{R}^{\varepsilon}_{\text{boundary}}$ are denoted by}
\begin{equation}
    \begin{aligned}
        \mathcal{R}^{\varepsilon}_{\text{residual}}   & = \; \frac{\lambda_1}{N_1^{(1)}} \sum_{i=1}^{N_1^{(1)}} | \varepsilon^2 \partial_t r^{\text{NN}}_{\theta}(t_i,x_i,v_i) + \varepsilon^2 v\partial_x j^{\text{NN}}_{\theta}(t_i,x_i,v_i)                                                                                     \\
                                                      & \quad \quad \quad \quad \quad \quad - (\rho^{\text{NN}}_{\theta}(t_i,x_i) - r^{\text{NN}}_{\theta}(t_i,x_i,v_i)) |^2                                                                                                                                                       \\
                                                      & + \frac{\lambda_2}{N_1^{(2)}} \sum_{i=1}^{N_1^{(2)}}  |\varepsilon^2 \partial_t j^{\text{NN}}_{\theta}(t_i,x_i,v_i) + v \partial_x r^{\text{NN}}_{\theta}(t_i,x_i,v_i) - (-j^{\text{NN}}_{\theta}(t_i,x_i,v_i)) |^2                                                        \\
                                                      & + \frac{\lambda_3}{N_1^{(3)}}\sum_{i=1}^{N_1^{(3)}} | \partial_t \rho^{\text{NN}}_{\theta}(t_i,x_i) +  \left \langle v\partial_x j^{\text{NN}}_{\theta} \right \rangle(t_i,x_i) |^2,                                                                                       \\
        \mathcal{R}^{\varepsilon}_{\text{constraint}} & = \frac{\lambda_4}{N_2}\sum_{i=1}^{N_2} |\rho^{\text{NN}}_{\theta}(t_i,x_i) - \left \langle r^{\text{NN}}_{\theta} \right \rangle(t_i,x_i) |^2,                                                                                                                            \\
        \mathcal{R}^{\varepsilon}_{\text{initial}}    & = \frac{\lambda_5}{N_3}\sum_{i=1}^{N_3}  |\rho^{\text{NN}}_{\theta}(0, x_i) - \left \langle f_{0} \right \rangle (x_i)|^2 + \frac{\lambda_6}{N_4} \sum_{i=1}^{N_4} |\mathcal{I}(r^{\text{NN}}_{\theta} + \varepsilon j^{\text{NN}}_{\theta})(x_i,v_i) - f_{0}(x_i,v_i)|^2, \\
        \mathcal{R}^{\varepsilon}_{\text{boundary}}   & = \frac{\lambda_7}{N_5} \sum_{i=1}^{N_5}  |\mathcal{B}(r^{\text{NN}}_{\theta} + \varepsilon j^{\text{NN}}_{\theta})(t_i,x_i,v_i) - F_{\text{B}}(t_i,x_i,v_i)|^2 .
    \end{aligned}
\end{equation}
Here, $N_1^{(1)}, N_1^{(2)}, N_1^{(3)}, N_2, N_3, N_4, N_5$ are the number of uniform sampling points of
domains $|\mathcal{T} \times \mathcal{D} \times \Omega|, |\mathcal{T} \times \mathcal{D} \times \Omega|, |\mathcal{T} \times \mathcal{D}|, |\mathcal{T} \times \mathcal{D}|, |\mathcal{D}|, |\mathcal{D} \times \Omega|, |\mathcal{T} \times\partial \mathcal{D} \times \Omega|$.

Case I and II are studied for problem 1 with constant scattering coefficient 1 for $\varepsilon = 10^{-3}, 10^{-8}$.
Furthermore, we focus on an uncertainty quantification problem for $\varepsilon = 10^{-5}$ in Case III.

    {\bf{Case I.} Inflow condition with $\varepsilon = 10^{-3}$}

Consider $\varepsilon=10^{-3}$, and the inflow boundary conditions are defined as $F_L(v) = 1$ and $F_R(v) = 0$, indicating the specified values of the distribution function at the left and right boundaries, respectively. The initial condition is set as $f_0(x, v) = 0$, denoting the initial distribution function.
Note that the function $f$ exhibits a discontinuity at $t = 0$ due to the conditions $F_L(v) = 1, F_R(v) = 0$, and $f_0(x, v) = 0$.
To enhance numerical performance, $\rho^{\text{NN}}_{\theta}$ can be further designed to inherently satisfy the initial condition.:
\begin{equation}
    \rho^{\text{NN}}_{\theta}(t, x) := t \cdot \exp \left( -\tilde{\rho}^{\text{NN}}_{\theta}(t, x)\right) \approx \rho(t, x).
\end{equation}

\begin{figure}[ht]
    \centering\includegraphics[width=0.45\textwidth]{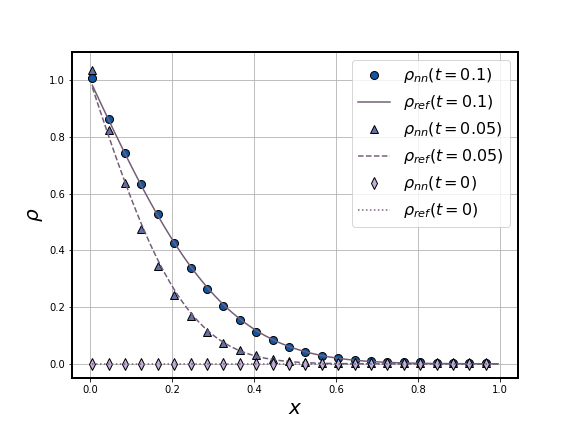}
    \caption{Problem 1---Case I. Plot of density $\rho$ at $t = 0, 0.05, 0.1$: APNNs (marker) vs. Ref (line).
        $\varepsilon = 10^{-3}$ and neural networks are FCNet with units $[2, 128, 128, 128, 128, 1]$ for $\rho$ and $[3, 256, 256, 256, 256, 1]$ both for $r$ and $j$.
        Batch size is $512$ in domain, $1024 \times 2$ on boundary and $512$ on initial condition.
        $\lambda_7 = 10$ and other $\lambda$'s are set to  be 1.
        Relative $\ell^2$ error of APNNs is $9.87 \times 10^{-3}$. }\label{fig:dirichlet10}
\end{figure}

Fig. \ref{fig:dirichlet10} depicts the estimated density, denoted as $\rho$, using APNNs in comparison to the reference solution at time instances $t = 0, 0.05, 0.1$.
It is evident that the approximated solutions exhibit favorable accuracy at both $t = 0, 0.05$ and $t = 0.1$.

The numerical performance of enforcing the initial condition and the soft constraint $\rho = \left \langle r \right \rangle$ is deliberated as follows.
The ensuing phenomena have been observed in all the experiments we have conducted.

\begin{figure}[ht]
    \centering\includegraphics[width=0.45\textwidth]{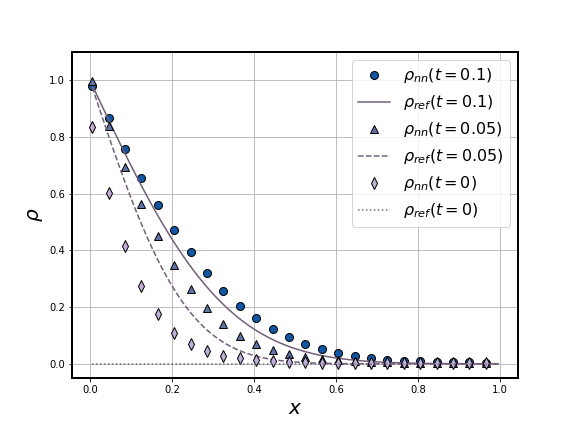}
    \caption{Problem 1---Case I. Plot of density $\rho$ at $t = 0, 0.05, 0.1$: APNNs (marker) vs. Ref (line).
        $\varepsilon = 10^{-3}$. And the units of neural networks are $[2, 128, 128, 128, 128, 1]$ for $\rho$ and $[3, 256, 256, 256, 256, 1]$ both for $r$ and $j$. }\label{fig:enforce_ic}
\end{figure}

Fig. \ref{fig:enforce_ic} illustrates the performance resulting from the improper enforcement of the initial condition.
Owing to the inadequate approximation of the initial condition and the influence of the boundary layer, erroneous solutions are obtained at time instances $t = 0, 0.05, 0.1$.

\begin{figure}[ht]
    \centering\includegraphics[width=0.45\textwidth]{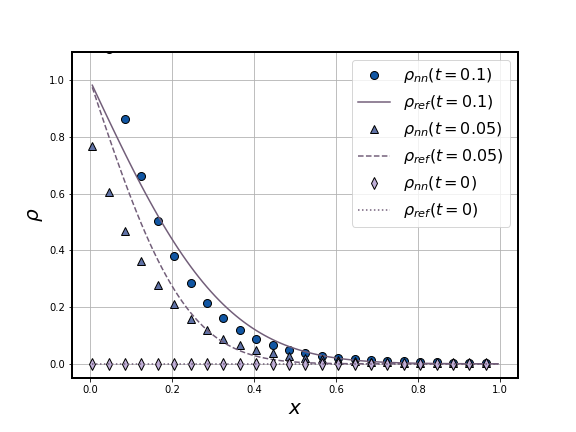}
    \caption{Problem 1---Case I. Plot of density $\rho$ at $t = 0, 0.05, 0.1$: APNNs (marker) vs. Ref (line).
        $\varepsilon = 10^{-3}$.
        And the units of neural networks are $[2, 128, 128, 128, 128, 1]$ for $\rho$ and $[3, 256, 256, 256, 256, 1]$ both for $r$ and $j$.}\label{fig:soft_constraint}
\end{figure}

Fig. \ref{fig:soft_constraint} illustrates the performance in the absence of  the constraint equation $\rho = \left \langle r \right \rangle$ in the loss function, while ensuring the exact satisfaction of the initial condition. However, it is evident that the solutions at time $t = 0.05, 0.1$ are incorrect. Consequently, we incorporate this constraint into our APNN loss function.

{\bf{Case II.} The homogeneous Dirichlet boundary condition with $\varepsilon = 10^{-8}$}

Let $\varepsilon = 10^{-8}$,
boundary condition be $F_L(v) = F_R(v) = 0$,
and initial condtion be
\begin{equation}
    f_0(x, v) = g(x) \cdot \frac{3}{\sqrt{2\pi}}e^{-\frac{{(3v)}^2}{2}},
\end{equation}
where
\begin{equation}
    g(x) = 1 + \sin \left ( 2 \pi x - \frac{\pi}{2} \right ).
\end{equation}

\begin{figure}[ht]
    \centering\includegraphics[width=0.45\textwidth]{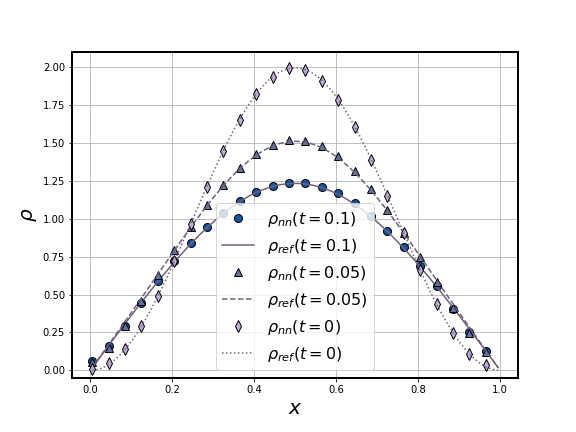}
    \caption{Problem 1---Case II.\@ Plot of density $\rho$ at $t = 0.05, 0.1$: APNNs (marker) vs. Ref (line).
        $\varepsilon = 10^{-8}$ and the units of neural networks are $[2, 128, 128, 128, 128, 128, 1]$ for $\rho$ and $[3, 256, 256, 256, 256, 256, 1]$ both for $r$ and $j$.
        Batch size is $1024$ in domain, $512 \times 2$ on boundary and $512$ on initial condition.
        $\lambda_3 = \lambda_4 = \lambda_7 = 10$ and others are set to  be 1.
        Relative $\ell^2$ error of APNNs is $1.25 \times 10^{-2}$. }\label{fig:dirichlet00}
\end{figure}

Fig. \ref{fig:dirichlet00} depicts the estimated density, denoted as $\rho$, obtained through the utilization of APNNs alongside the reference solution at three distinct time instances, namely $t = 0, 0.05, 0.1$.
In this particular scenario, it is noteworthy that the function $f$ exhibits a continuous behavior without any discontinuous jumps, while possessing a non-constant initial value.
Evidently, this methodology has proven to yield commendable outcomes, as illustrated in Fig. \ref{fig:dirichlet00}.

{\bf{Case III.} Uncertainty quantification problem with inflow condition and $(\varepsilon = 10^{-5})$}

Next, we contemplate an uncertainty quantification problem ($\varepsilon = 10^{-5}$):
\begin{equation}
    \varepsilon \partial_t f + v\partial_x f = \frac{\sigma_S(\bm z)}{\varepsilon}\left ( \frac{1}{2} \int_{-1}^1 f \, dv' - f \right ), \quad x_L < x < x_R, \quad -1 \leq v \leq 1,
\end{equation}
with scattering coefficient
{
\begin{equation}
    \sigma_S(\bm z) = 1 + \frac{1}{10} \sum_{i=1}^{10} \prod_{j=1}^{i} z^j, \;\bm{z} = (z^1, z^2, \ldots, z^{10}) \sim \mathcal{U}({[-1, 1]}^{10}),
\end{equation}}
and in-flow boundary condition $F_L(v) = 1, F_R(v) = 0$, and initial condition $f_0(x, v) = 0$.

In order to address this issue, we incorporate the {10}-dimensional stochastic vector ${\bm z}$ as a constituent of the input components of the deep neural networks for $\rho, r, j$.
Similarly, it is possible to derive the empirical APNN risk for the even-odd system of linear transport equation with uncertainties in the following manner
\begin{equation}
    \mathcal{R}^{\varepsilon}_{\text{APNN, uq}} = \mathcal{R}^{\varepsilon}_{\text{residual}} + \mathcal{R}^{\varepsilon}_{\text{constraint}} + \mathcal{R}^{\varepsilon}_{\text{initial}} + \mathcal{R}^{\varepsilon}_{\text{boundary}},
\end{equation}
{where $\mathcal{R}^{\varepsilon}_{\text{residual}}, \mathcal{R}^{\varepsilon}_{\text{constraint}}, \mathcal{R}^{\varepsilon}_{\text{initial}},  \mathcal{R}^{\varepsilon}_{\text{boundary}}$ are denoted by
\begin{equation*}
    \begin{aligned}
        \mathcal{R}^{\varepsilon}_{\text{residual}} = & \; \frac{\lambda_1}{N_1^{(1)}} \sum_{i=1}^{N_1^{(1)}} | \varepsilon^2 \partial_t r^{\text{NN}}_{\theta}({\bm z_i}, t_i,x_i,v_i) + \varepsilon^2 v\partial_x j^{\text{NN}}_{\theta}({\bm z_i}, t_i,x_i,v_i)                                                                                                              \\
                                                      & \quad \quad \quad \quad \quad - \sigma_S(\bm z_i)(\rho^{\text{NN}}_{\theta}({\bm z_i}, t_i,x_i) - r^{\text{NN}}_{\theta}({\bm z_i}, t_i,x_i,v_i)) |^2                                                                                                                                                                   \\
                                                      & \; + \frac{\lambda_2}{N_1^{(2)}} \sum_{i=1}^{N_1^{(2)}}  |\varepsilon^2 \partial_t j^{\text{NN}}_{\theta}({\bm z_i}, t_i,x_i,v_i) + v \partial_x r^{\text{NN}}_{\theta}({\bm z_i}, t_i,x_i,v_i) - (-\sigma_S(\bm z_i) j^{\text{NN}}_{\theta}({\bm z_i}, t_i,x_i,v_i)) |^2                                               \\
                                                      & \; + \frac{\lambda_3}{N_1^{(3)}}\sum_{i=1}^{N_1^{(3)}} | \partial_t \rho^{\text{NN}}_{\theta}({\bm z_i}, t_i,x_i) +  \left \langle v\partial_x j^{\text{NN}}_{\theta} \right \rangle({\bm z_i}, t_i,x_i) |^2,   \\
        \mathcal{R}^{\varepsilon}_{\text{constraint}} & = \frac{\lambda_4}{N_2}\sum_{i=1}^{N_2} |\rho^{\text{NN}}_{\theta}({\bm z_i}, t_i,x_i) - \left \langle r^{\text{NN}}_{\theta} \right \rangle({\bm z_i}, t_i,x_i) |^2,                                                                                                                            \\
        \mathcal{R}^{\varepsilon}_{\text{initial}} =  & \; \frac{\lambda_5}{N_3}\sum_{i=1}^{N_3}  |\rho^{\text{NN}}_{\theta}({\bm z_i}, 0, x_i) - \left \langle f_{0} \right \rangle ({\bm z_i}, x_i)|^2 + \frac{\lambda_6}{N_4} \sum_{i=1}^{N_4} |\mathcal{I}(r^{\text{NN}}_{\theta} + \varepsilon j^{\text{NN}}_{\theta})({\bm z_i}, x_i,v_i) - f_{0}({\bm z_i}, x_i,v_i)|^2, \\
        \mathcal{R}^{\varepsilon}_{\text{boundary}} = & \; \frac{\lambda_7}{N_5} \sum_{i=1}^{N_5}  |\mathcal{B}(r^{\text{NN}}_{\theta} + \varepsilon j^{\text{NN}}_{\theta})({\bm z_i}, t_i,x_i,v_i) - F_{\text{B}}({\bm z_i}, t_i,x_i,v_i)|^2 .
    \end{aligned}
\end{equation*}
Here, $N_1^{(1)}, N_1^{(2)}, N_1^{(3)}, N_2, \cdots, N_5$ is the number of sample points of corresponding domain.}

{
To assess the numerical performance, we evaluate the value of $\rho$ at specific time instances, namely $t = 0.05, 0.1$.
This evaluation is based on computing the expected outcome over $10^3$ simulation iterations, considering ${\bm z}$ as a vector comprising elements $z^1, z^2, \ldots, z^{10}$ with each $z^i$ generated uniformly from $[-1, 1]$.}

\begin{figure}[htbp]
    \centering\includegraphics[width=0.5\textwidth]{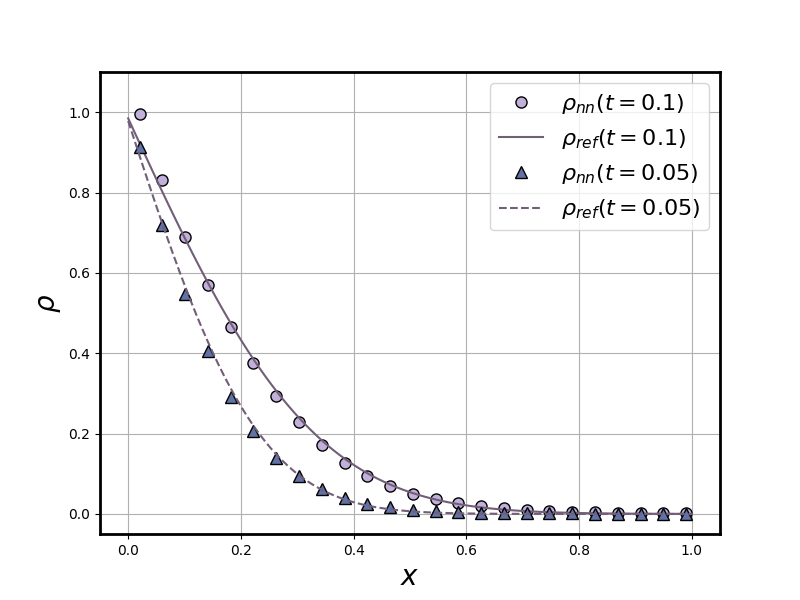}
    \caption{{Problem 1---Case III.\@ Plot of density $\rho$ by taking expectation for $\bm{z}$ at $t = 0.05, 0.1$ for APNNs (marker) and Ref (line).
        $\varepsilon = 10^{-5}, \sigma_S(\bm z) = 1 + \frac{1}{10} \sum_{i=1}^{10} \prod_{j=1}^{i} z^j$ and
        the units of neural networks are $[22, 64, 128, 256, 256, 128, 64, 1]$ for $\rho$ and $[23, 128, 256, 512, 512, 256, 128, 1]$ for $r, j$.
        Batch size is $1024$ in domain, $512 \times 2$ for boundary condition and $256$ for initial condition.
        $\lambda_5 = \lambda_7 = 10$ and others are set to be 1.
        Relative $\ell^2$ error of APNNs is $3.90 \times 10^{-2}$.} }\label{fig:uq}
\end{figure}

Fig. \ref{fig:uq} exhibits the anticipated density $\rho$, trained using APNNs, alongside the reference solutions at time $t = 0.05, 0.1$.
The results vividly demonstrate the exceptional capabilities of APNNs when dealing with high-dimensional problems.

\subsection{Problem 2: APNNs for solving Boltzmann-BGK equations}
Consider the Boltzmann-BGK equation
\begin{equation}
    \partial_t f + v \cdot \nabla_x f = \frac{1}{\varepsilon} \left ( M(U) - f \right ),  \quad v \in \mathbb{R}.
\end{equation}
and recall the system of BGK model for our APNN is
\begin{equation*}
    \left \{
    \begin{aligned}
         & \varepsilon \left ( \partial_t f + v \partial_x f \right ) = M(U) - f, \\
         & \partial_t U
        + \nabla_x \cdot \left \langle v m f \right \rangle = 0,                  \\
         & U = \left \langle m f \right \rangle .
    \end{aligned}
    \right.
\end{equation*}

Initial local thermodynamics equilibrium is considered for all tests:
\begin{equation*}
    f(0, x, v) = \frac{\rho_0}{{(2 \pi T_0)}^{\frac{1}{2}}} \exp \left ( - \frac{|v - u_0|^2}{2 T_0} \right ) := f_0 (x, v).
\end{equation*}

For these tests, a classical ``Riemann-like'' problem with Sod like initial data is considered.
\begin{equation*}
    \begin{aligned}
         & \rho_0(t, x) = \rho_L, & \quad u_0(t, x) = u_L, & \quad T_0(t, x) = T_L, & \quad x \leq 0, \\
         & \rho_0(t, x) = \rho_R, & \quad u_0(t, x) = u_R, & \quad T_0(t, x) = T_R, & \quad x > 0.
    \end{aligned}
\end{equation*}

By saying ``Riemann-like'', we allow the initial data to be a smoothed curve connecting two constants of each macroscopic quantities instead of just two constants (see the following test cases I, II, III, IV, V), i.e.,
\begin{equation*}
    \begin{aligned}
         & \rho_0(t, x_L) = \rho_L, \quad \rho_0(t, x_R) = \rho_R, \\
         & u_0(t, x_L) = u_L, \quad u_0(t, x_R) = u_R,             \\
         & T_0(t, x_L) = T_L, \quad T_0(t, x_R) = T_R,
    \end{aligned}
\end{equation*}
where
\begin{equation*}
    W_L = \left (
    \rho_L, \;
    u_L, \;
    T_L
    \right )^T, \;
    W_R = \left (
    \rho_R, \;
    u_R, \;
    T_R
    \right )^T.
\end{equation*}
Due to the finite propagation speed, the above conditions hold as long as the waves initially started from the Riemann-like initial data have not reached the computational boundaries in the duration of the computation time.

Numerically, one can construct four DNNs as follows:
\begin{equation}
    f^{\text{NN}}_{\theta}(t, x, v) := \ln \left(1 + \exp (\tilde{f}^{\text{NN}}_{\theta}(\bar{t}, x, \bar{v})) \right) > 0,
\end{equation}
and
\begin{equation}
    \begin{aligned}
         & \rho ^{\text{NN}}_{\theta}(t, x) := \rho_L^{\frac{x_R - x}{x_R - x_L}} \cdot \rho_R^{\frac{x - x_L}{x_R - x_L}} \cdot \exp \left ( (x-x_L)(x_R - x) \cdot  \tilde{\rho}^{\text{NN}}_{\theta}(\bar{t}, x)\right ) > 0, \\
         & u ^{\text{NN}}_{\theta}(t, x) := \sqrt{(x-x_L)(x_R - x)} \cdot \tilde{u}^{\text{NN}}_{\theta}(\bar{t}, x) + \frac{x_R - x}{x_R - x_L} u_L + \frac{x - x_L}{x_R - x_L} u_R,                                            \\
         & T ^{\text{NN}}_{\theta}(t, x) := T_L^{\frac{x_R - x}{x_R - x_L}} \cdot T_R^{\frac{x - x_L}{x_R - x_L}} \cdot \exp \left ( (x-x_L)(x_R - x) \cdot  \tilde{T}^{\text{NN}}_{\theta}(\bar{t}, x)\right ) > 0,
    \end{aligned}
\end{equation}
which $\rho ^{\text{NN}}_{\theta}, u ^{\text{NN}}_{\theta}, T ^{\text{NN}}_{\theta}$ automatically satisfy the boundary constraint which we find can improve the numerical performance.
In this problem, to keep $f$ positive, $\ln (1 + \exp(\cdot))$ is applied for constructing $f^{\text{NN}}_{\theta}$.
The advantage of this structure lies in the parity of magnitudes between $f^{\text{NN}}_{\theta}$ and $\tilde{f}^{\text{NN}}_{\theta}$ when their magnitude exceeds a certain threshold.

The APNN empirical risk for the system of Boltzmann-BGK equation is
\begin{equation}
{
    \mathcal{R}^{\varepsilon}_{\text{APNN, BGK}} = \mathcal{R}^{\varepsilon}_{\text{residual}} + 
    \mathcal{R}^{\varepsilon}_{\text{claw}} +
    \mathcal{R}^{\varepsilon}_{\text{constraint}} + \mathcal{R}^{\varepsilon}_{\text{boundary}} + \mathcal{R}^{\varepsilon}_{\text{initial}},}
\end{equation}
{where $\mathcal{R}^{\varepsilon}_{\text{residual}}, \mathcal{R}^{\varepsilon}_{\text{claw}}, \mathcal{R}^{\varepsilon}_{\text{constraint}}, \mathcal{R}^{\varepsilon}_{\text{boundary}}, \mathcal{R}^{\varepsilon}_{\text{initial}}$ are denoted by
\begin{small}
    \begin{equation}
        \begin{aligned}
            \mathcal{R}^{\varepsilon}_{\text{residual}}   & = \frac{\lambda_1}{N_1} \sum_{i=1}^{N_1} | \varepsilon (\partial_t f^{\text{NN}}_{\theta}(t_i,x_i,v_i) +  v_i \nabla_x f^{\text{NN}}_{\theta}(t_i,x_i,v_i)) - \left ( M(U^{\text{NN}}_{\theta}) - f^{\text{NN}}_{\theta} \right )(t_i,x_i,v_i) |^2, \\
            \mathcal{R}^{\varepsilon}_{\text{claw}}       & = \frac{{\lambda_2}}{N_2} \sum_{i=1}^{N_2} |\partial_t U^{\text{NN}}_{\theta}(t_i,x_i) + \nabla_x \left \langle v m f^{\text{NN}}_{\theta} \right \rangle(t_i,x_i) |^2 ,                                                                         \\
            \mathcal{R}^{\varepsilon}_{\text{constraint}} & = \frac{{\lambda_3}}{N_3} \sum_{i=1}^{N_3}  |U^{\text{NN}}_{\theta}(t_i,x_i) - \left \langle m f^{\text{NN}}_{\theta} \right \rangle(t_i,x_i) |^2,                                                                                               \\        
            \mathcal{R}^{\varepsilon}_{\text{boundary}}   & = \frac{{\lambda_4}}{N_4} \sum_{i=1}^{N_4} |W^{\text{NN}}_{\theta}(t_i, x_L) - W_L|^2 + |W^{\text{NN}}_{\theta}(t_i, x_R) - W_R|^2, \\
            \mathcal{R}^{\varepsilon}_{\text{initial}}    & = \frac{{\lambda_5}}{N_5} \sum_{i=1}^{N_5} |W^{\text{NN}}_{\theta}(0, x_i) - W_0(x_i)|^2 + \frac{\lambda_{6}}{N_6} \sum_{i=1}^{N_6} |f^{\text{NN}}_{\theta}(0, x_i, v_i) - f_0(x_i, v_i)|^2. 
        \end{aligned}
    \end{equation}
\end{small}}
Here, $N_i (i = 1,\cdots,6)$ is the number of sample points of corresponding domain.

These cases are studied in Problem 2.

    {\bf{Case I:}} $\varepsilon=10^{-3}$

\begin{equation*}
    \begin{aligned}
         & \rho_0(x) = 1.5 + (0.625 - 1.5) \cdot \frac{\sin(\pi x) + 1}{2}, \\
         & u_0(x) =0,                                                       \\
         & T_0(x) = 1.5 + (0.75 - 1.5) \cdot \frac{\sin(\pi x) + 1}{2}.
    \end{aligned}
\end{equation*}

\begin{figure}[ht]
    \subfigure[{The integrals of approximate $f_\text{nn}$ vs.\ reference solutions}. {\it Left}: $t = 0$ and {\it Right}: $t = 0.1$.]
    {
        \includegraphics[width=0.45\textwidth]{ 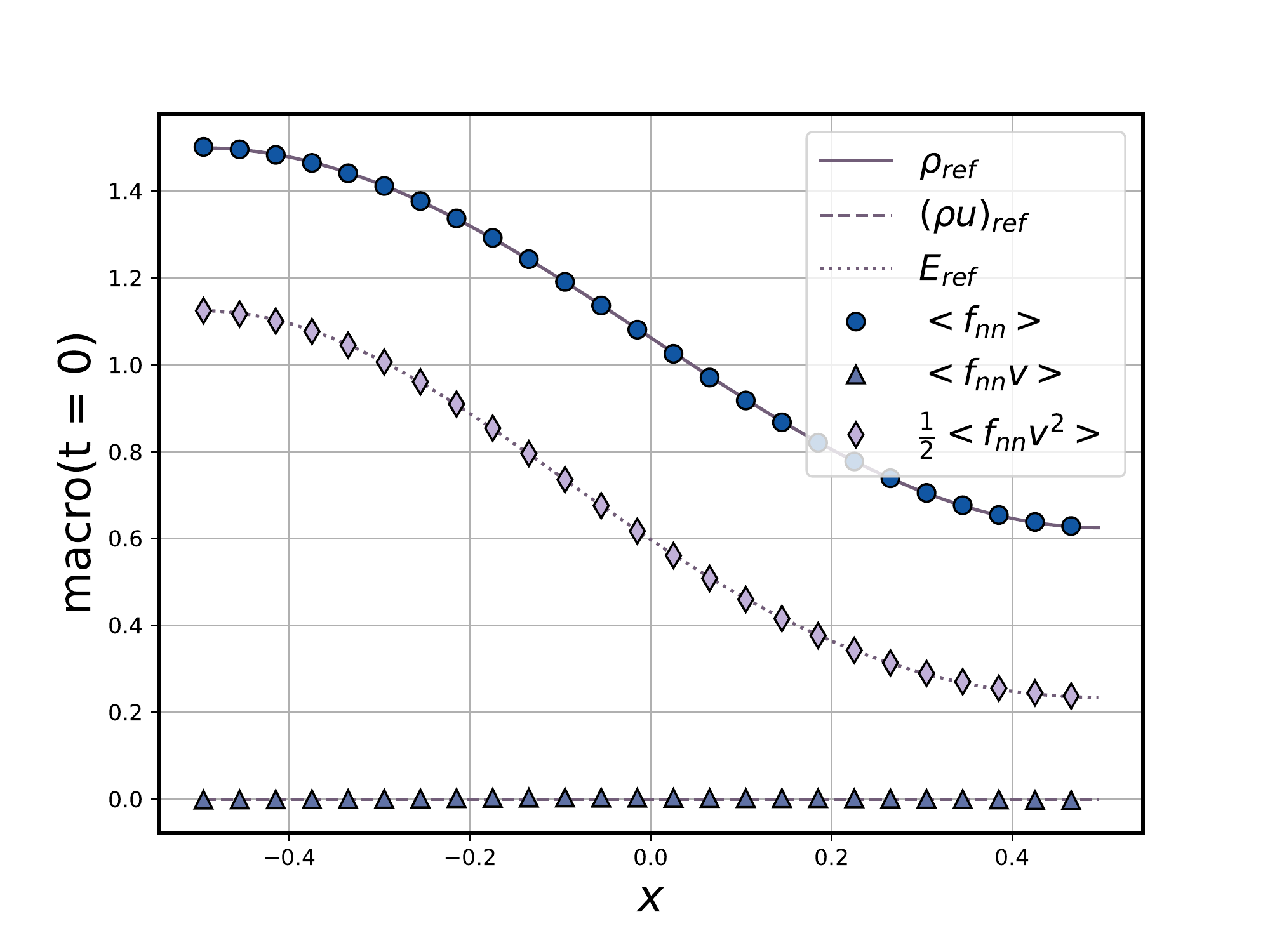}
        \includegraphics[width=0.45\textwidth]{ 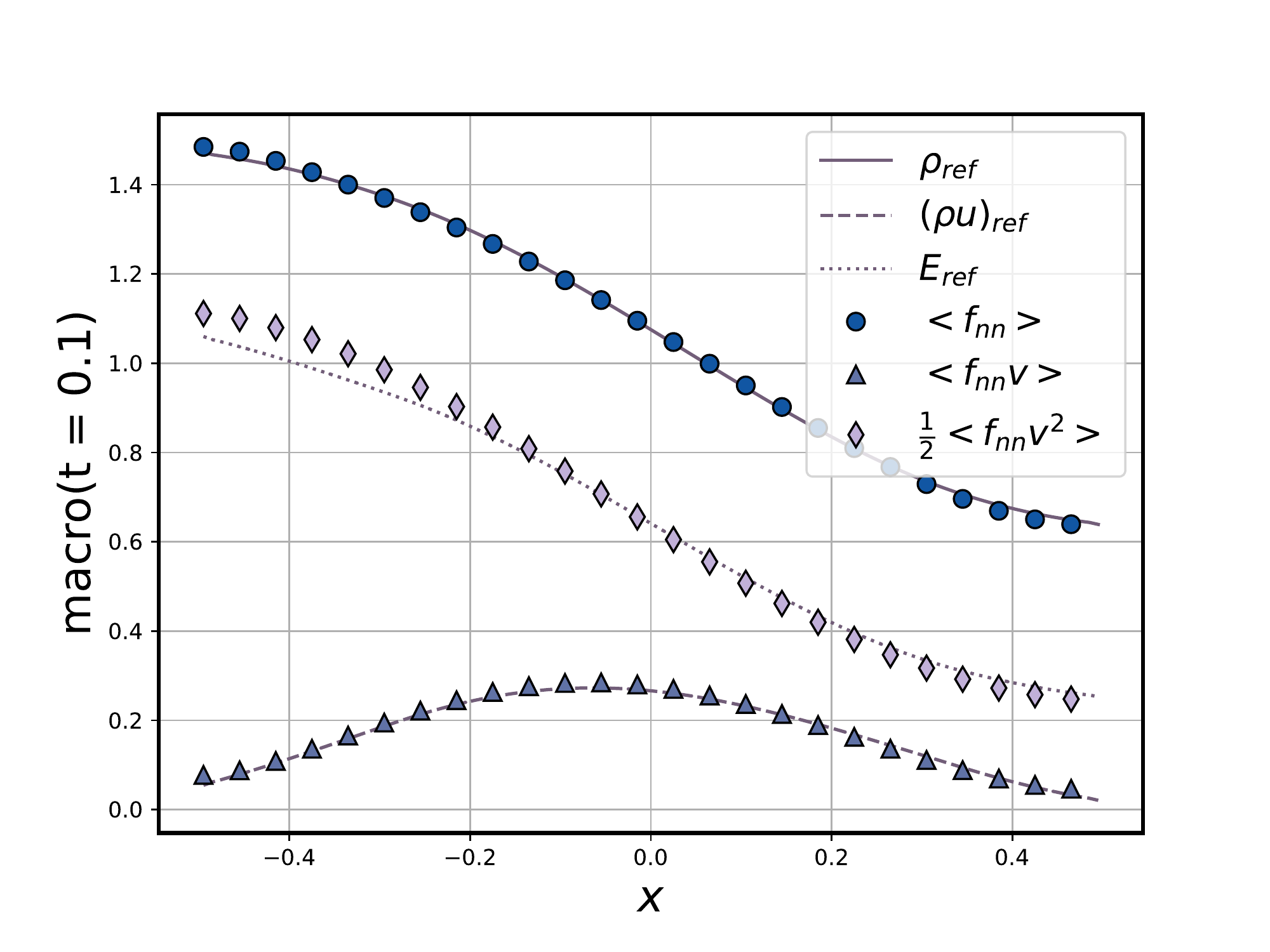}
    }

    \subfigure[{The approximate $\rho_\text{nn}, u_\text{nn}, T_\text{nn}$} vs.\ reference solutions. {\it Left}: $t = 0$ and {\it Right}: $t = 0.1$.]
    {
        \includegraphics[width=0.45\textwidth]{ 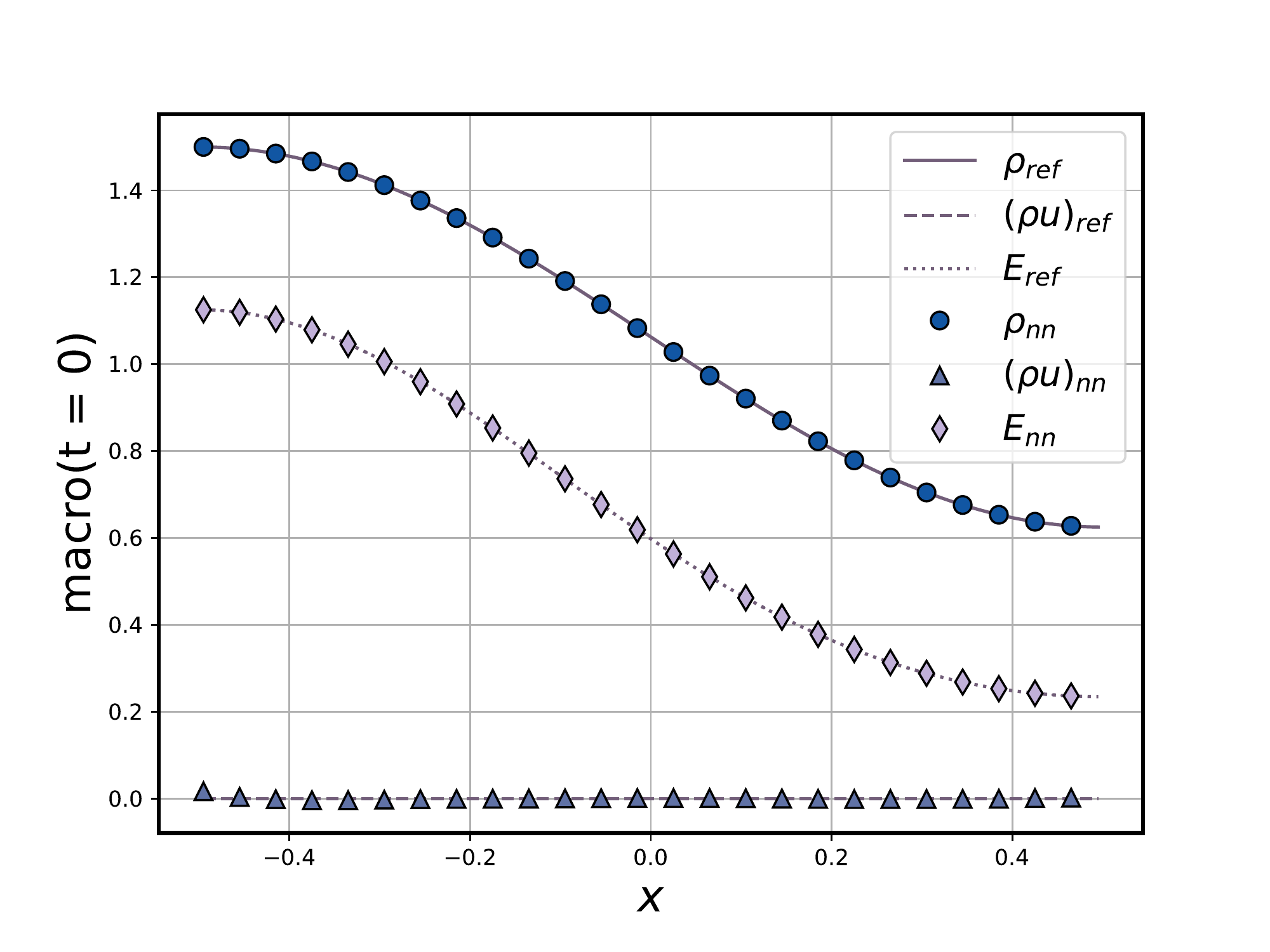}
        \includegraphics[width=0.45\textwidth]{ 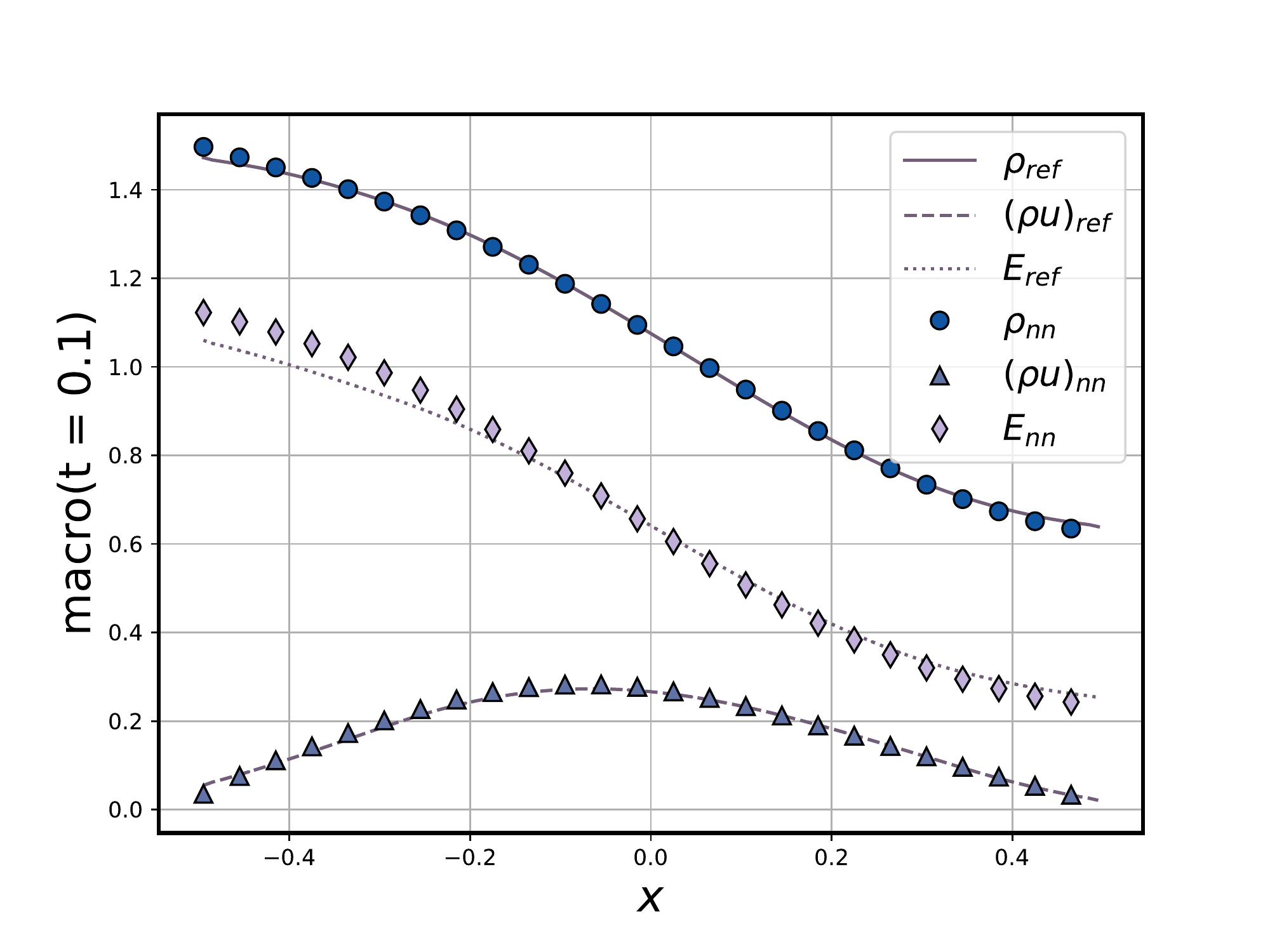}
    }
    \caption{Problem 2---Case I. Plot of density, momentum and energy at time $t = 0, 0.1$: Approximated by APNNs (marker) vs. Ref (line).
        $\varepsilon = 10^{-3}$ and the units of neural networks are $[3, 128, 128, 128, 128, 128, 128, 1]$ for $f$ and $[2, 64, 64, 64, 64, 64, 64, 1]$ both for $\rho, u$ and $T$.
        Batch size is $512$ in domain, and $256$ on initial condition. ${{\lambda_5}} = (1, 10, 10), \lambda_{6} = 10$ and others are set to be $1$.
        For $t = 0$: mean square error of density, momentum and energy are $7.16\text{e-8}, 4.50\text{e-6}, 8.13\text{e-7}$.
        For $t = 0.1$: relative $l^2$ error of density, momentum and energy are $5.43\text{e-3}, 6.35\text{e-3}, 4.47\text{e-2}$.}\label{fig:example-sin-1}
\end{figure}

Fig. \ref{fig:example-sin-1} illustrates the graph representing approximate macroscopic properties at time $t = 0, 0.1$. 
It is evident that the approximate density, momentum, and energy exhibit superior performance compared to those derived from the approximate function $f$.

{\bf{Case II:}} $\varepsilon=10^{-3}$

\begin{equation*}
    \begin{aligned}
         & \rho_0(x) = 1.5 + (0.625 - 1.5) \cdot
        \frac{
            {\sin}^{\frac{3}{7}} (\pi x) + 1
        }
        {2},                                     \\
         & u_0(x) =0,                            \\
         & T_0(x) = 1.5 + (0.75 - 1.5) \cdot
        \frac{
            {\sin}^{\frac{3}{7}} (\pi x) + 1
        }
        {2}.
    \end{aligned}
\end{equation*}

\begin{figure}[ht]
    \subfigure[{The integrals of approximate $f_\text{nn}$ vs.\ reference solutions}. {\it Left}: $t = 0$ and {\it Right}: $t = 0.1$.]
    {
        \includegraphics[width=0.45\textwidth]{ 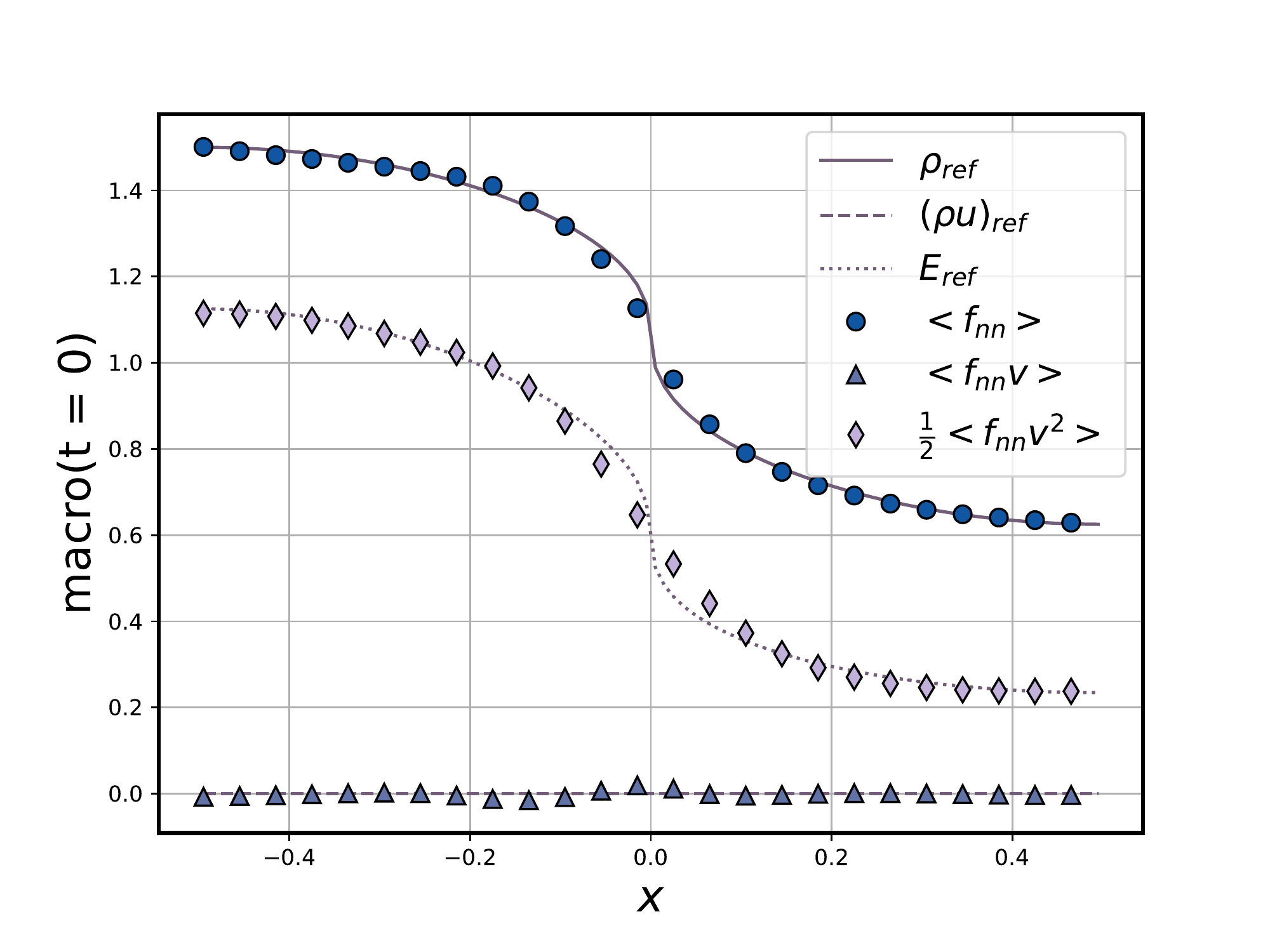}
        \includegraphics[width=0.45\textwidth]{ 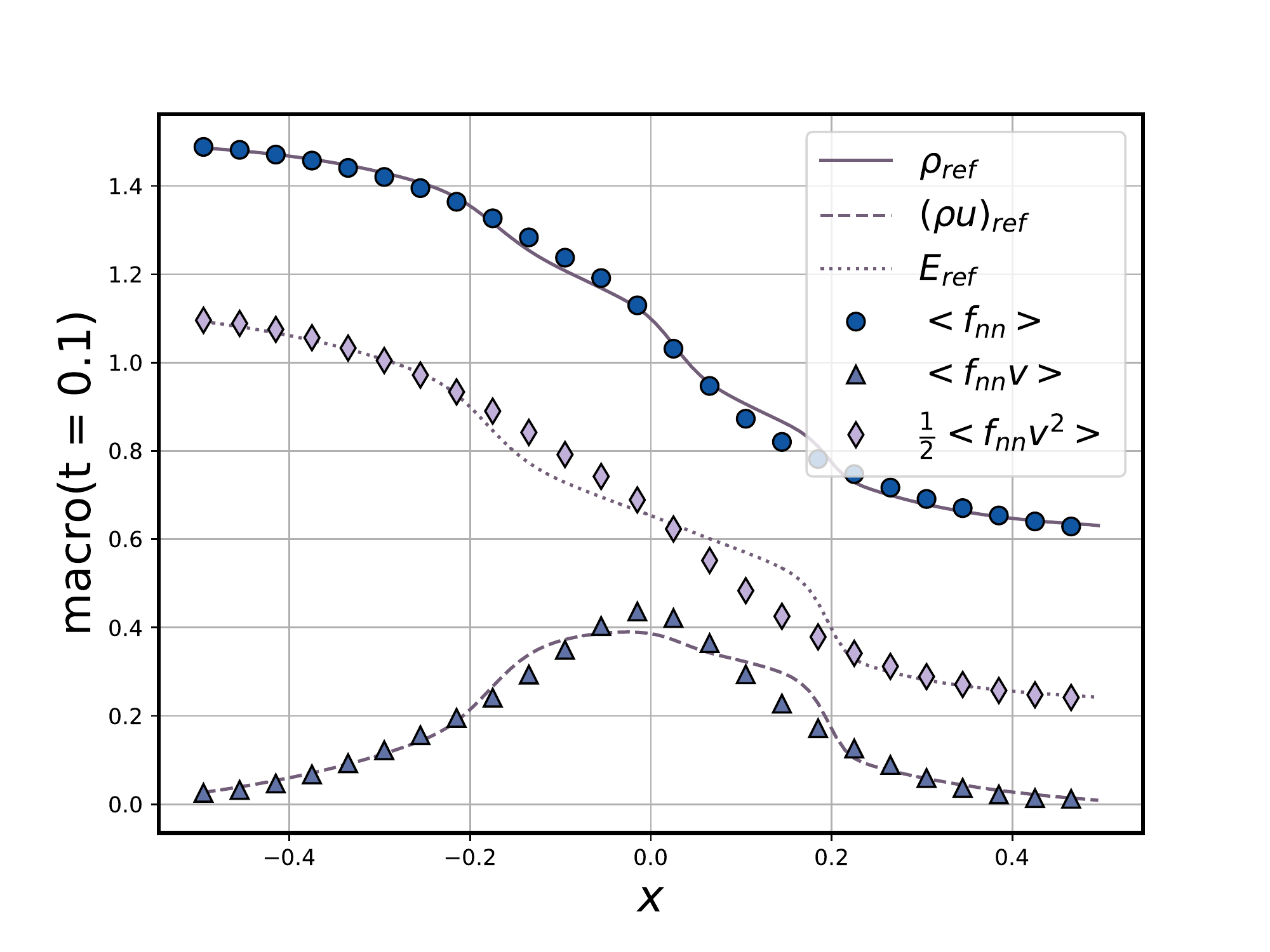}
    }

    \subfigure[{The approximate $\rho_\text{nn}, u_\text{nn}, T_\text{nn}$} vs.\ reference solutions. {\it Left: $t = 0$} and {\it Right: $t = 0.1$}.]
    {
        \includegraphics[width=0.45\textwidth]{ 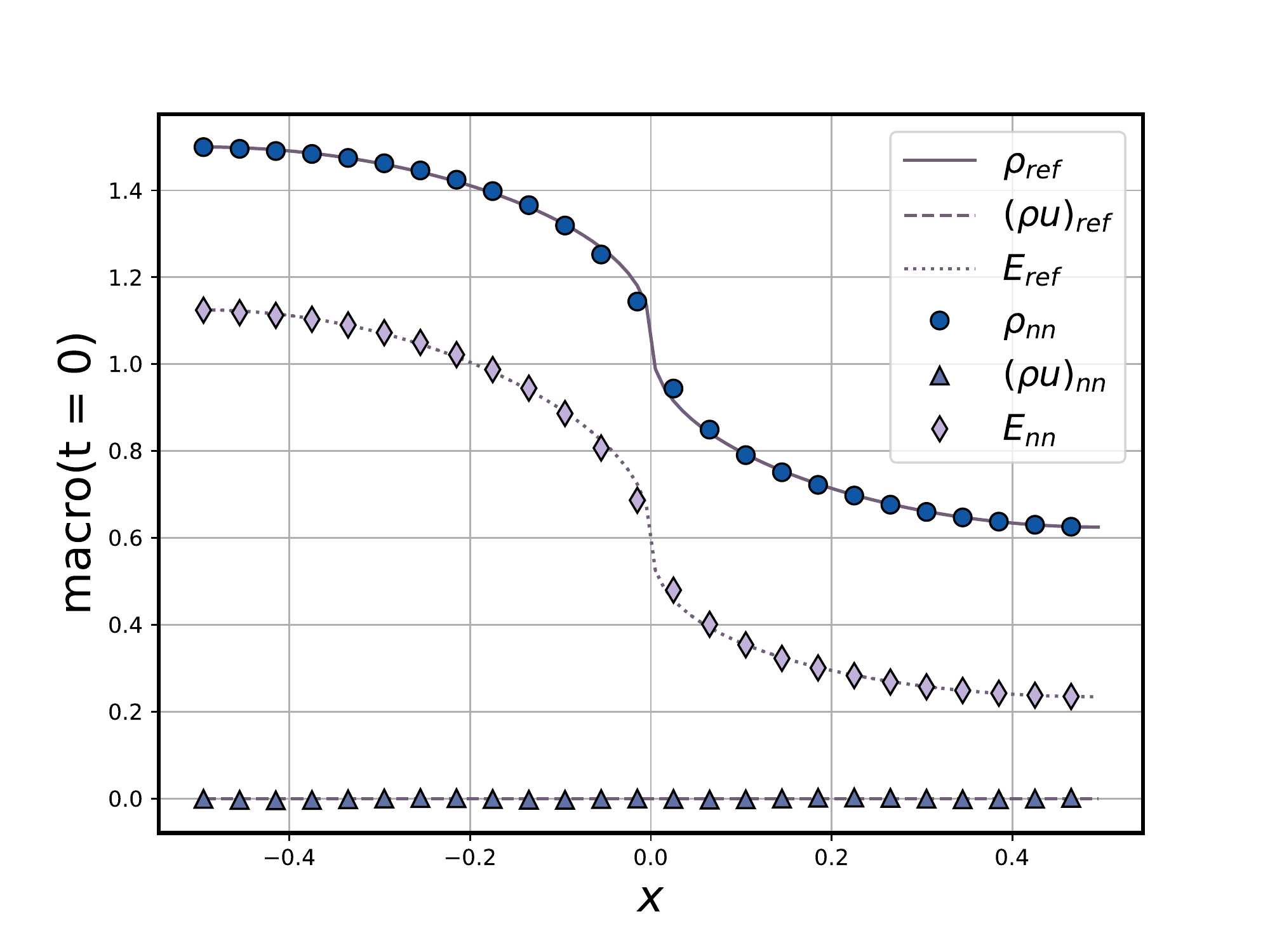}
        \includegraphics[width=0.45\textwidth]{ 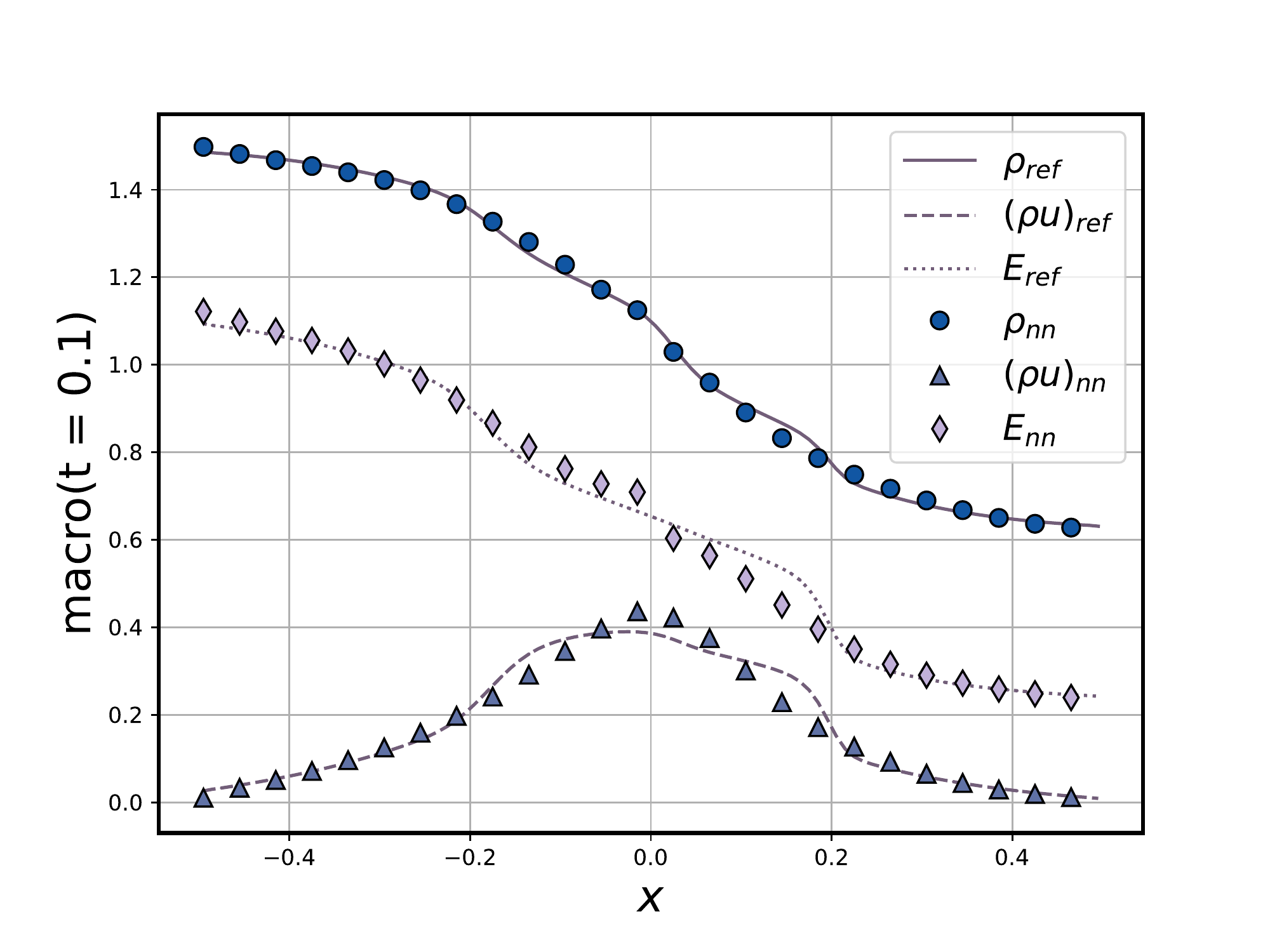}
    }
    \caption{Problem 2---Case II.\@ Plot of density, momentum and energy at time $t = 0, 0.1$: Approximated by APNNs (marker) vs. Ref (line).
        $\varepsilon = 10^{-3}$ and the units of neural networks are $[3, 128, 128, 128, 128, 128, 128, 1]$ for $f$ and $[2, 64, 64, 64, 64, 64, 64, 1]$ both for $\rho, u$ and $T$.
        Batch size is $512$ in domain, and $256$ on initial condition. ${{\lambda_5}} = (1, 10, 10)$ and others are set to be $1$.
        For $t = 0$: mean square error of density, momentum and energy are $1.29\text{e-4}, 6.34\text{e-6}, 4.41\text{e-5}$.
        For $t = 0.1$: relative $l^2$ error of density, momentum and energy are $1.36\text{e-2}, 2.00\text{e-2}, 3.99\text{e-2}$.}\label{fig:example-sin-2}
\end{figure}

Fig. \ref{fig:example-sin-2} illustrates the plot depicting the approximated macroscopic quantities at $t = 0, 0.1$. It is noteworthy that the steep slope of $\rho(0, x)$ and $T(0, x)$ at $x = 0$ poses a considerable challenge for obtaining accurate solutions through approximation. The key observation in this scenario is that approximating the moments of $f$, namely, $\rho, u, T$, is comparatively simpler than approximating $f$ itself. This assertion holds true across all the test cases.

{\bf{Case III:}} $\varepsilon=10^{-3}$

\begin{equation*}
    \begin{aligned}
         & \rho_0(x) = 1.5 + (0.625 - 1.5) \cdot
        \frac{
            \tanh (10 x) + 1
        }
        {2},                                     \\
         & u_0(x) =0,                            \\
         & T_0(x) = 1.5 + (0.75 - 1.5) \cdot
        \frac{
            \tanh (10 x) + 1
        }
        {2}.
    \end{aligned}
\end{equation*}

\begin{figure}[ht]
    \subfigure[{The integrals of approximate $f_\text{nn}$ vs.\ reference solutions}. {\it Left: $t = 0$} and {\it Right: $t = 0.1$}.]
    {
        \includegraphics[width=0.45\textwidth]{ 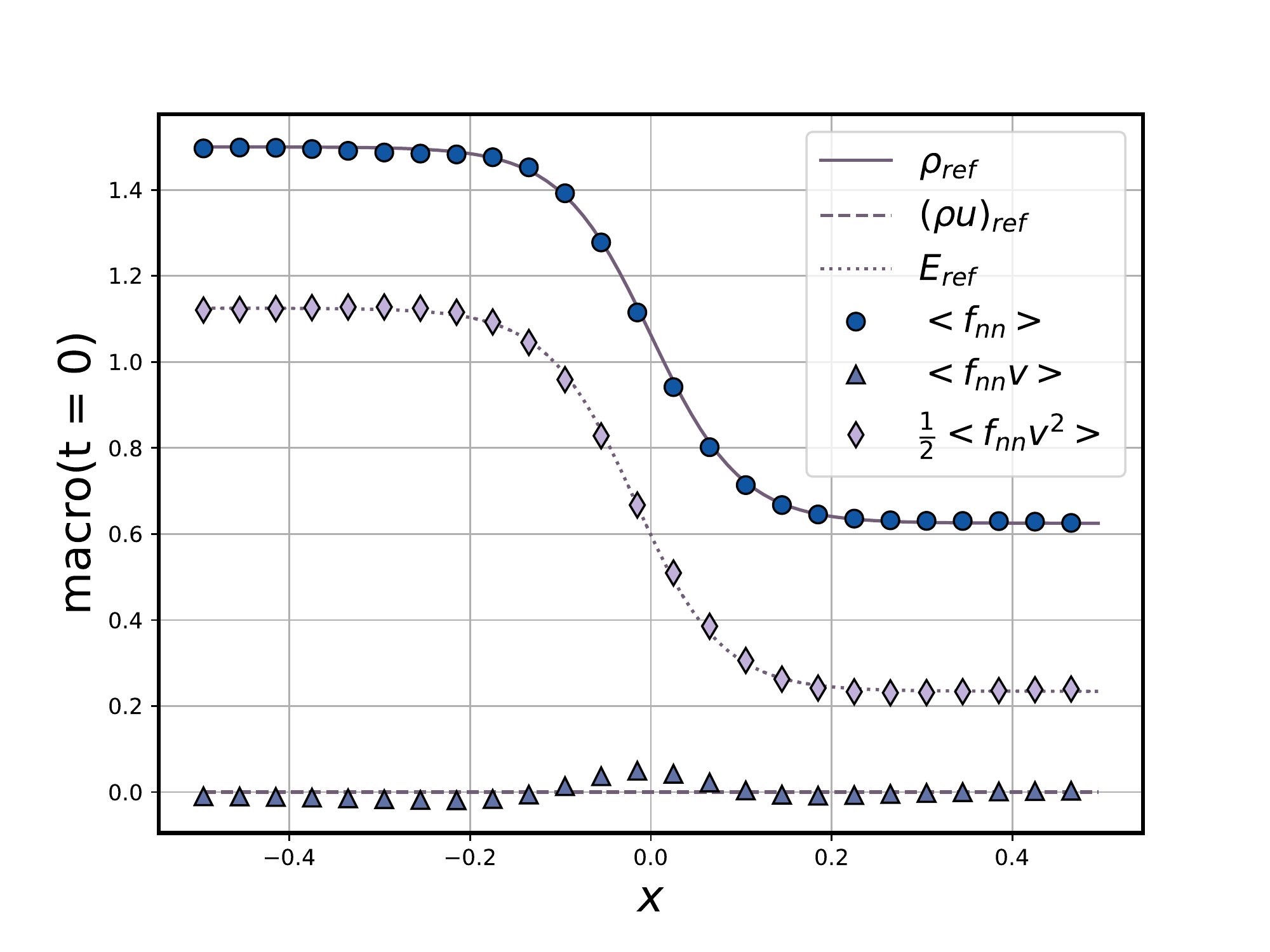}
        \includegraphics[width=0.45\textwidth]{ 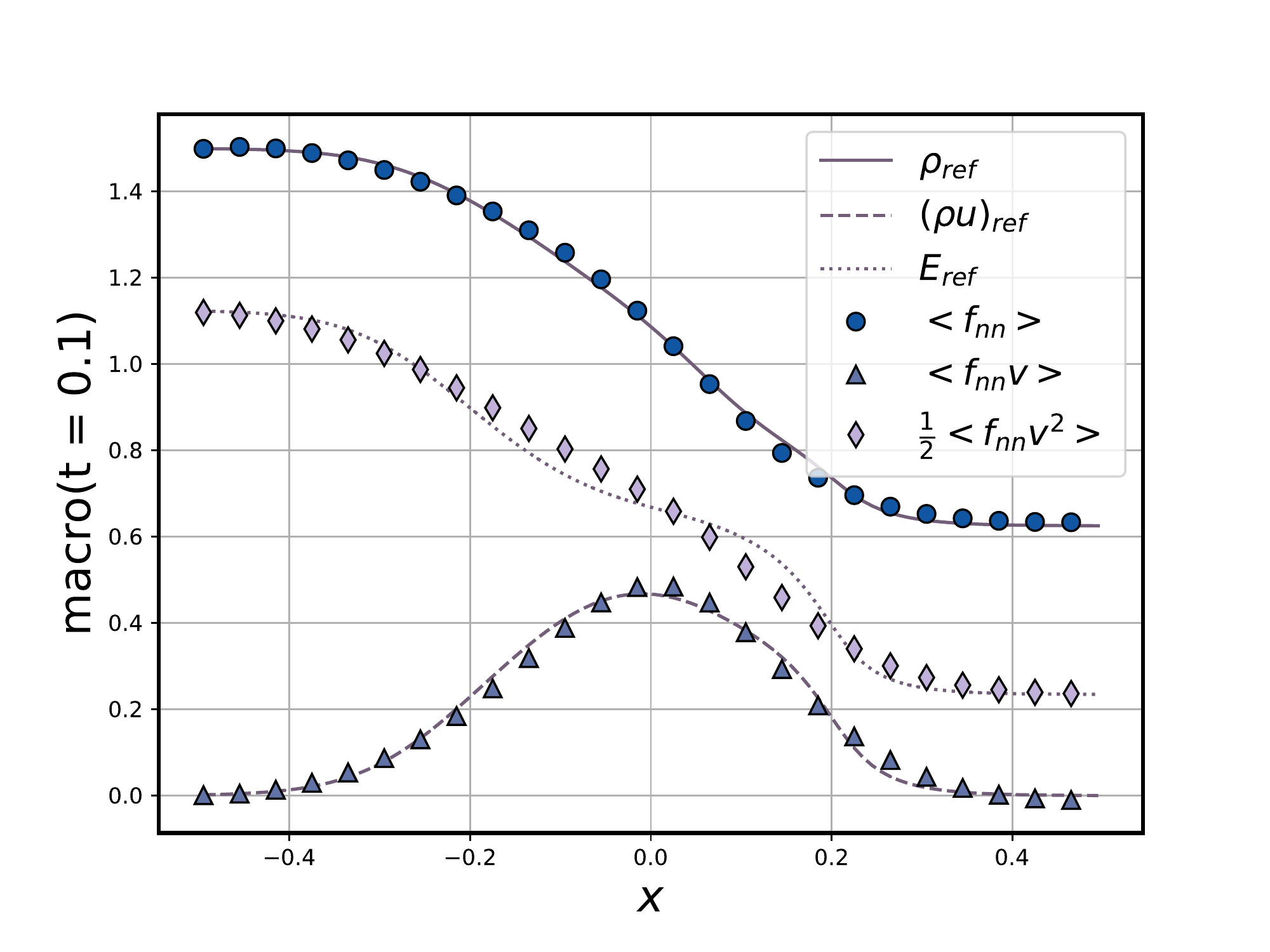}
    }

    \subfigure[{The approximate $\rho_\text{nn}, u_\text{nn}, T_\text{nn}$} vs.\ reference solutions. {\it Left: $t = 0$} and {\it Right: $t = 0.1$}.]
    {
        \includegraphics[width=0.45\textwidth]{ 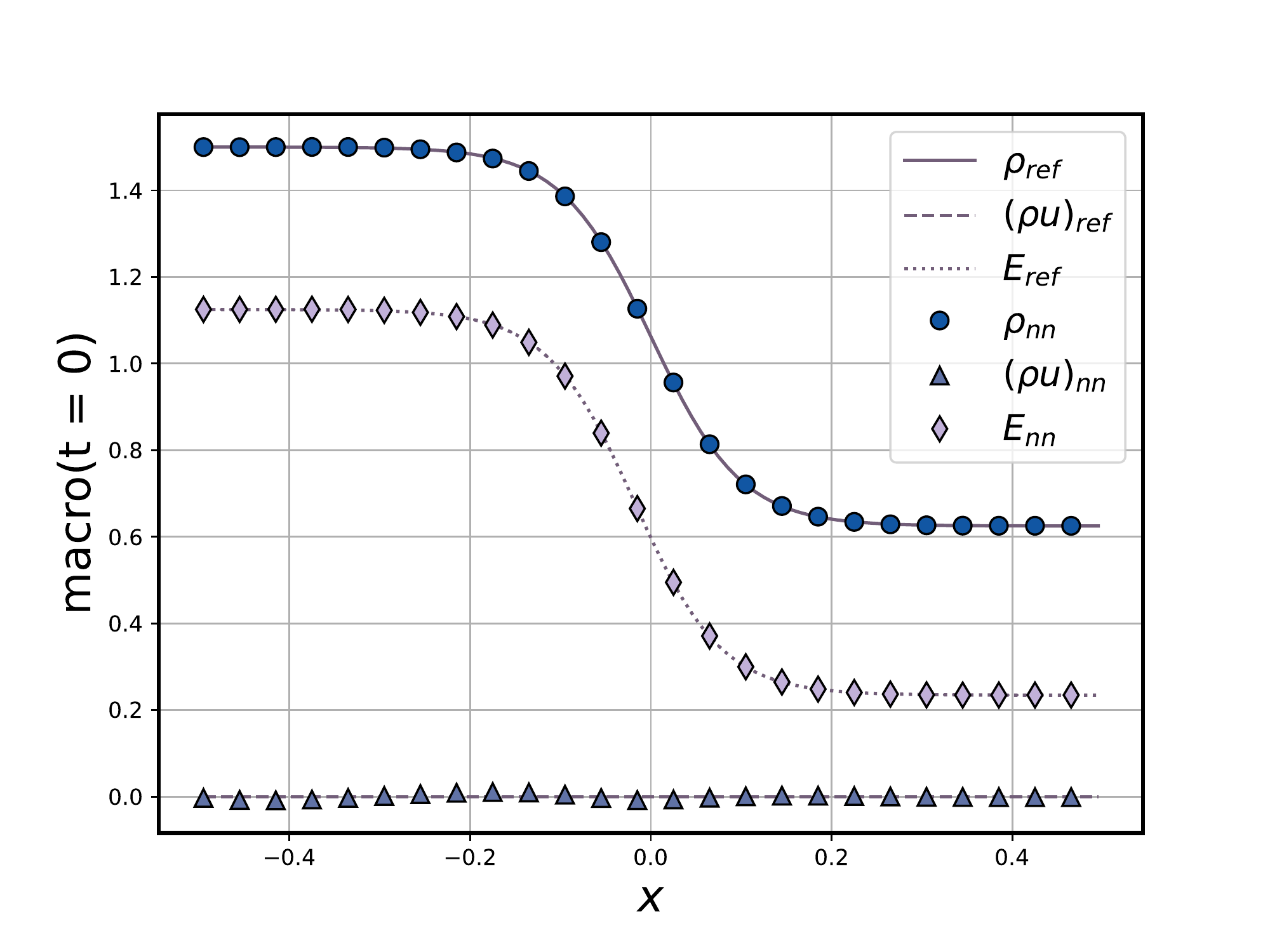}
        \includegraphics[width=0.45\textwidth]{ 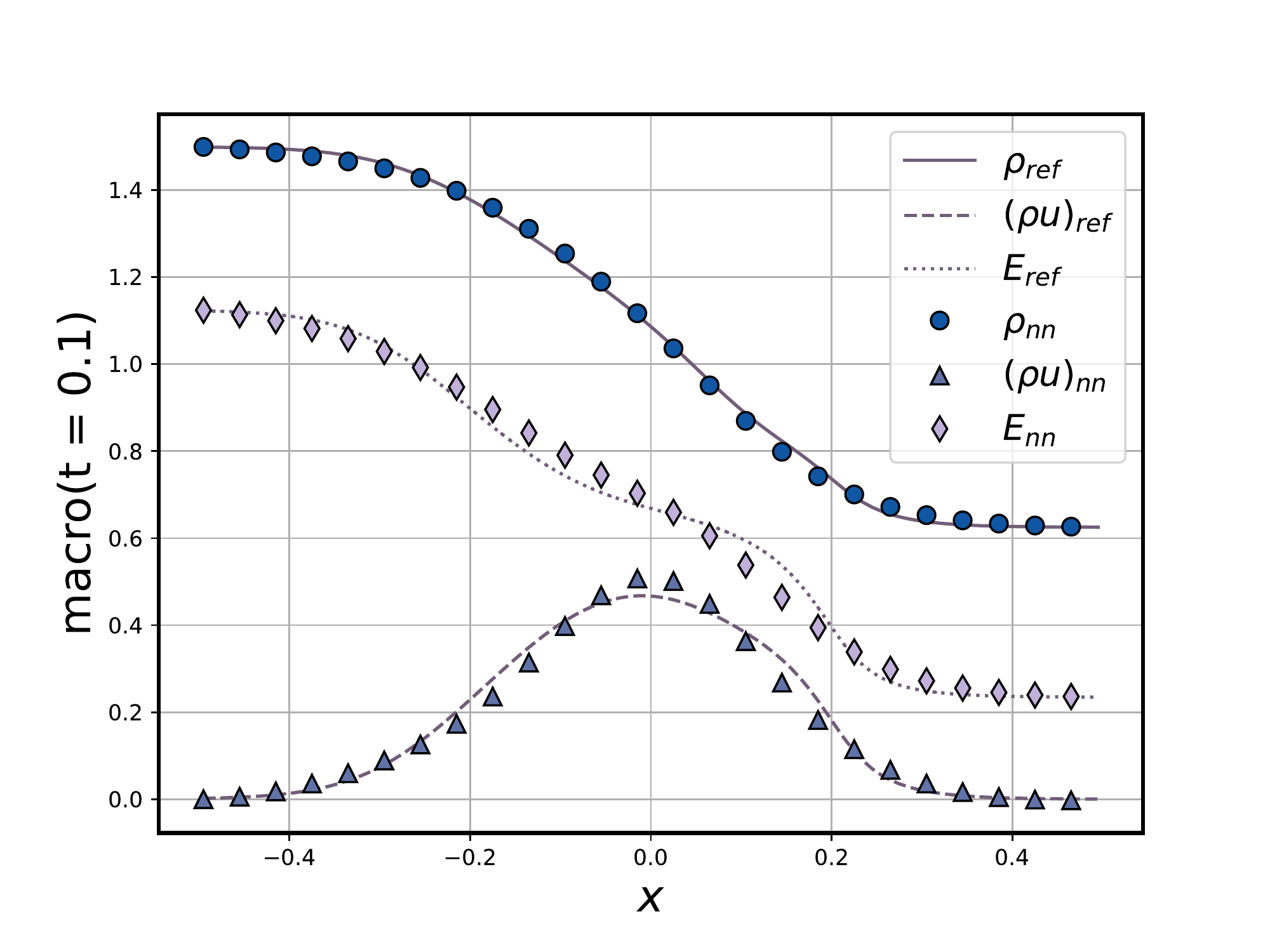}
    }
    \caption{Problem 2---Case III.\@ Plot of density, momentum and energy at time $t = 0, 0.1$: Approximated by APNNs (marker) vs. Ref (line).
        $\varepsilon = 10^{-3}$ and the units of neural networks are $[3, 128, 128, 128, 128, 128, 128, 1]$ for $f$ and $[2, 64, 64, 64, 64, 64, 64, 1]$ both for $\rho, u$ and $T$.
        Batch size is $512$ in domain, and $256$ on initial condition. ${{\lambda_5}} = (1, 10, 10), \lambda_6 = 10$ and others are set to be $1$.
        For $t = 0$: mean square error of density, momentum and energy are $4.87\text{e-8}, 1.22\text{e-6}, 3.29\text{e-8}$.
        For $t = 0.1$: relative $l^2$ error of density, momentum and energy are $6.19\text{e-3}, 1.78\text{e-2}, 3.60\text{e-2}$.}\label{fig:example-tanh-1}
\end{figure}

Fig. \ref{fig:example-tanh-1} depicts the graph illustrating the plot of approximate macroscopic quantities at $t = 0, 0.1$.
In this particular scenario, the gradient of $\rho(0, x)$ and $T(0, x)$ at $x = 0$ exhibits a seamless behavior.
It is evident that the approximate solutions align well with the reference solutions; however, the first moment computed by the approximate function $f$ yields unsatisfactory outcomes in the vicinity of $x = 0$.
This observation concurs with the findings presented in cases I and II.

{\bf{Case IV:}} $\varepsilon=1$

\begin{equation*}
    \begin{aligned}
         & \rho_0(x) = 1.5 + (0.625 - 1.5) \cdot
        \frac{
            \tanh (20 x) + 1
        }
        {2},                                     \\
         & u_0(x) =0,                            \\
         & T_0(x) = 1.5 + (0.75 - 1.5) \cdot
        \frac{
            \tanh (20 x) + 1
        }
        {2}.
    \end{aligned}
\end{equation*}

\begin{figure}[ht]
    \subfigure[{The integrals of approximate $f_\text{nn}$ vs.\ reference solutions}. {\it Left: $t = 0$} and {\it Right: $t = 0.1$}.]
    {
        \includegraphics[width=0.45\textwidth]{ 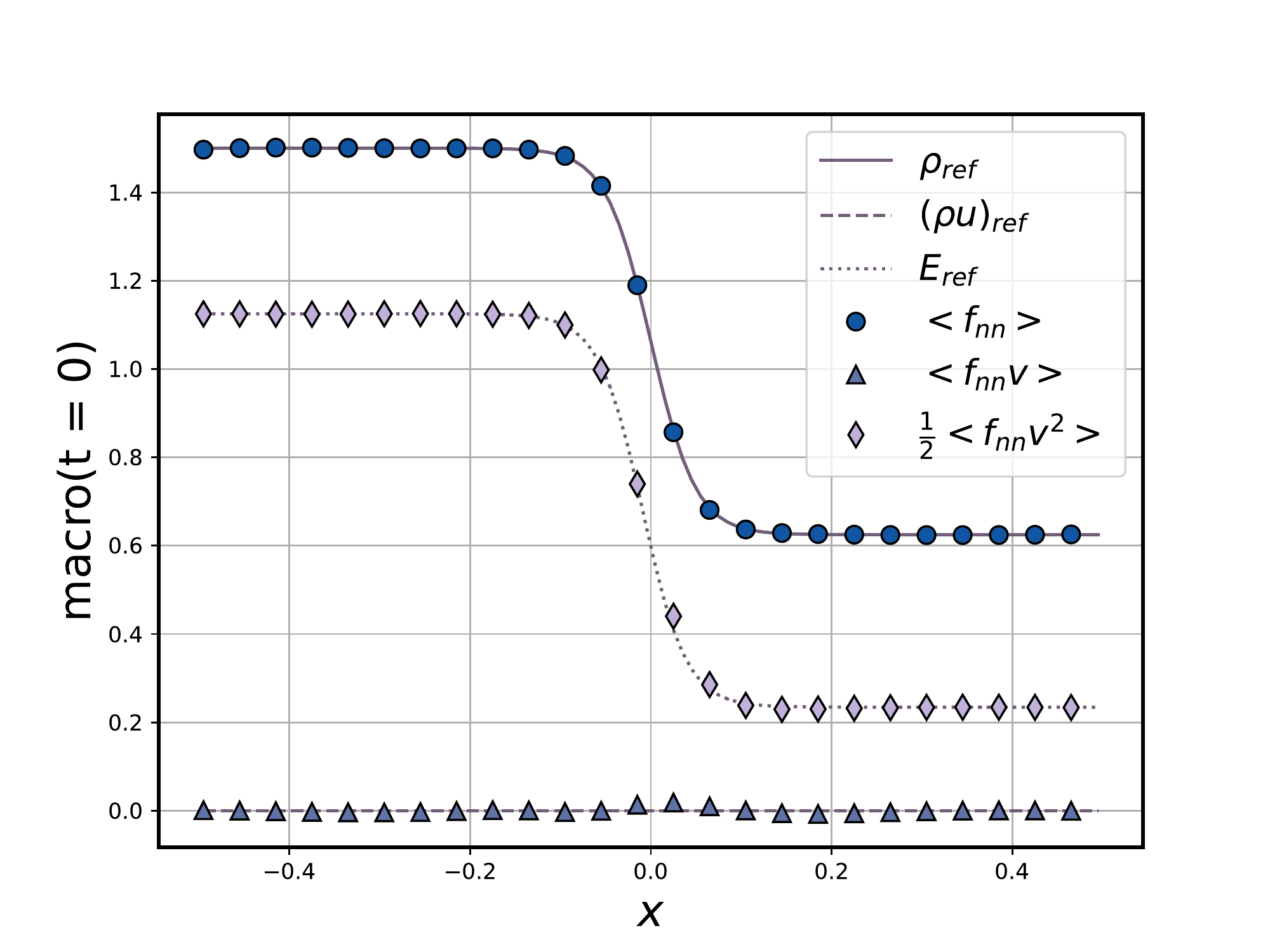}
        \includegraphics[width=0.45\textwidth]{ 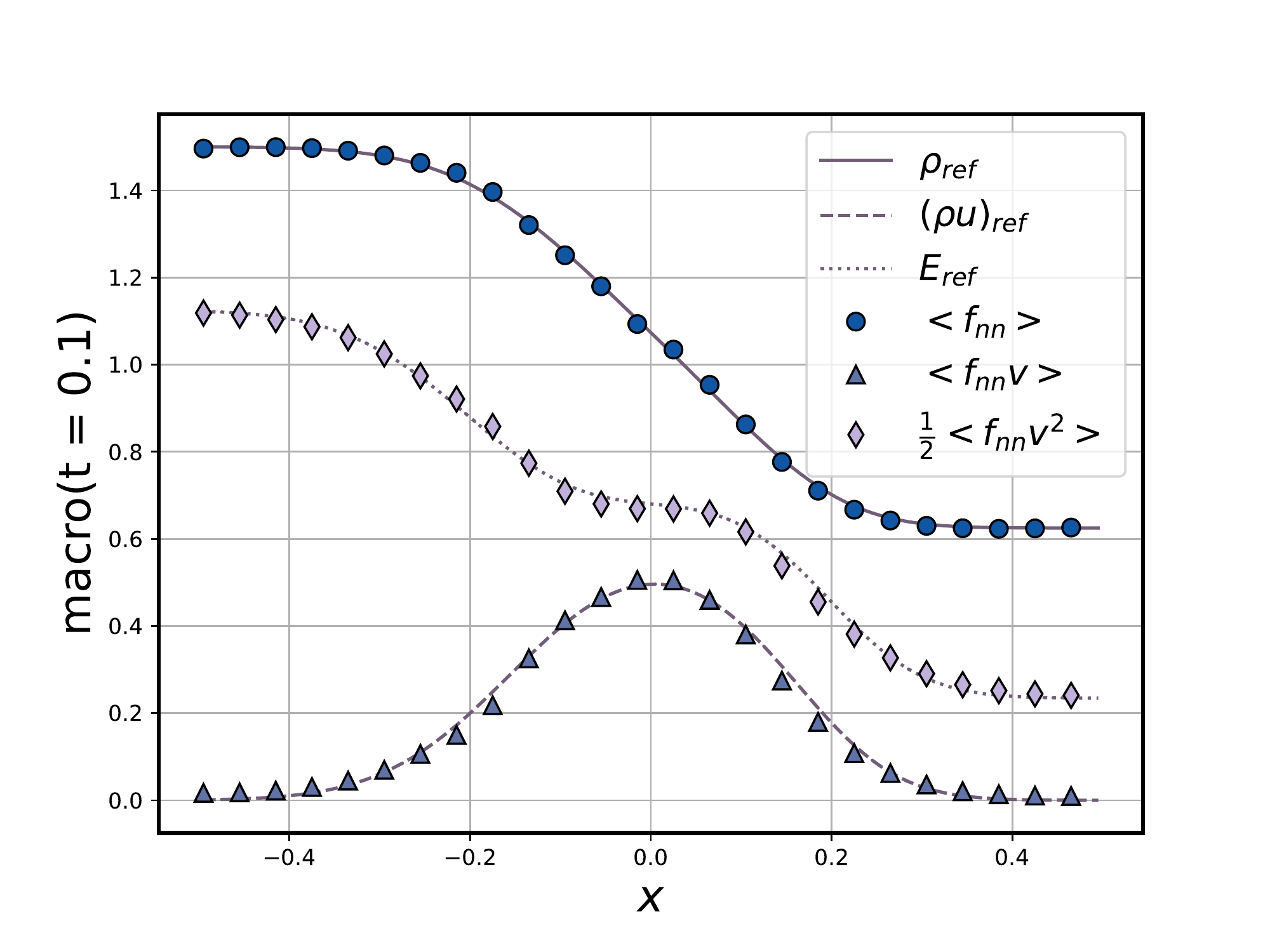}
    }
    \subfigure[{The approximate $\rho_\text{nn}, u_\text{nn}, T_\text{nn}$} vs.\ reference solutions. {\it Left: $t = 0$} and {\it Right: $t = 0.1$}.]
    {
        \includegraphics[width=0.45\textwidth]{ 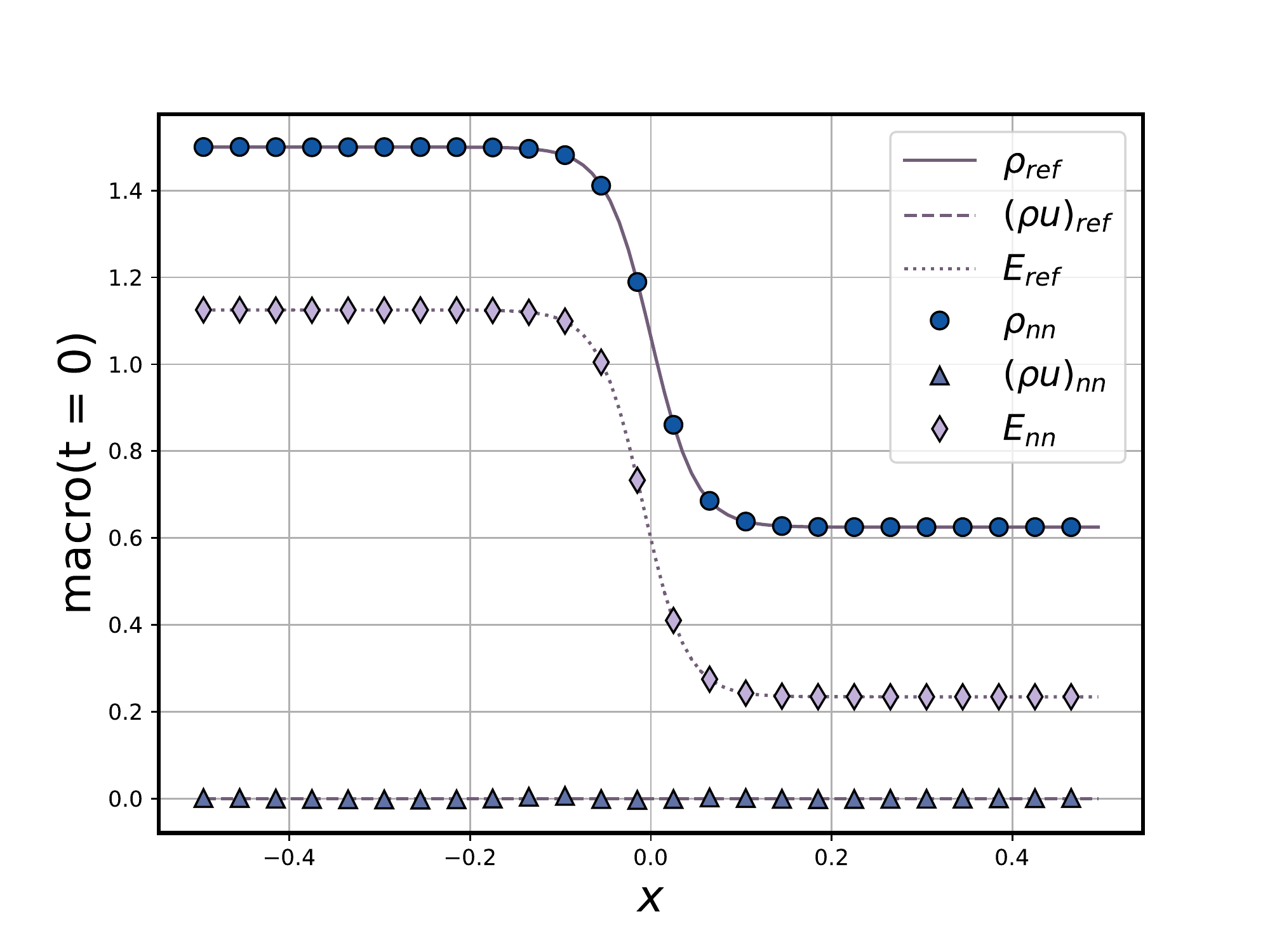}
        \includegraphics[width=0.45\textwidth]{ 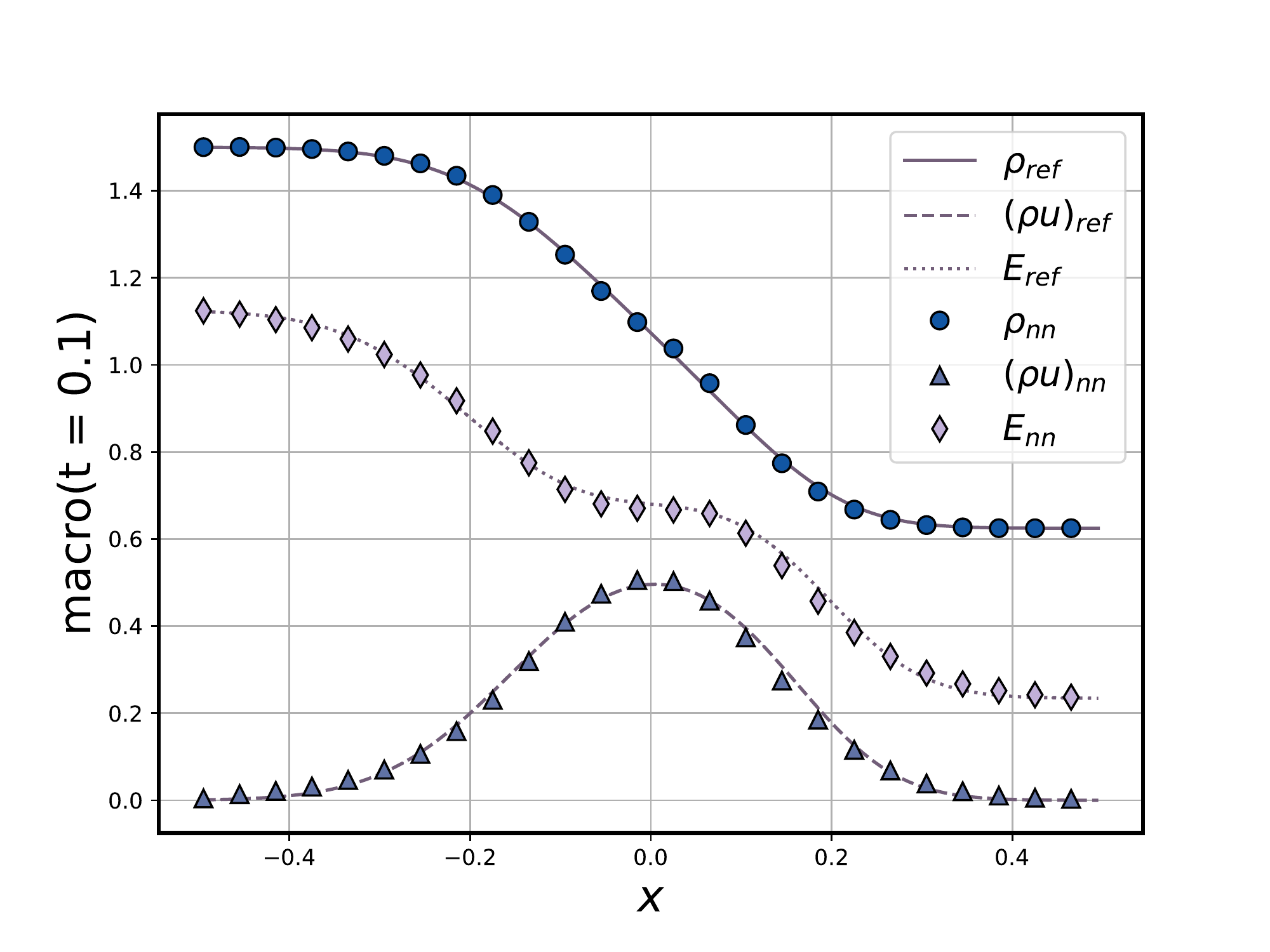}
    }
    \caption{Problem 2---Case IV.\@ Plot of density, momentum and energy at time $t = 0, 0.1$: Approximated by APNNs (marker) vs. Ref (line).
        $\varepsilon = 1$ and the units of neural networks are $[3, 128, 128, 128, 128, 128, 128, 1]$ for $f$ and $[2, 64, 64, 64, 64, 64, 64, 1]$ both for $\rho, u$ and $T$.
        Batch size is $512$ in domain, and $256$ on initial condition.
        ${{\lambda_5}} = (1, 10, 10), \lambda_6 = 10$ and others are set to be $1$.
        For $t = 0$: mean square error of density, momentum and energy are $9.89\text{e-8}, 2.34\text{e-6}, 2.08\text{e-7}$.
        For $t = 0.1$: relative $l^2$ error of density, momentum and energy are $6.41\text{e-3}, 8.72\text{e-3}, 1.62\text{e-2}$.}\label{fig:example-tanh-2}
\end{figure}

The plot displayed in Fig. \ref{fig:example-tanh-2} illustrates the representation of estimated macroscopic quantities at time instances $t = 0, 0.1$.
It is worth noting that this particular scenario closely resembles a realistic problem.
In this context, our proposed APNN method demonstrates commendable efficacy by offering favorable solutions when compared to the reference solutions.

\begin{figure}[ht]
    \subfigure[{{\it Left:} Case I with $\varepsilon=10^{-3}$} and {{\it Right:} Case II with $\varepsilon=10^{-3}$}.]
    {
    \includegraphics[width=0.45\textwidth]{ 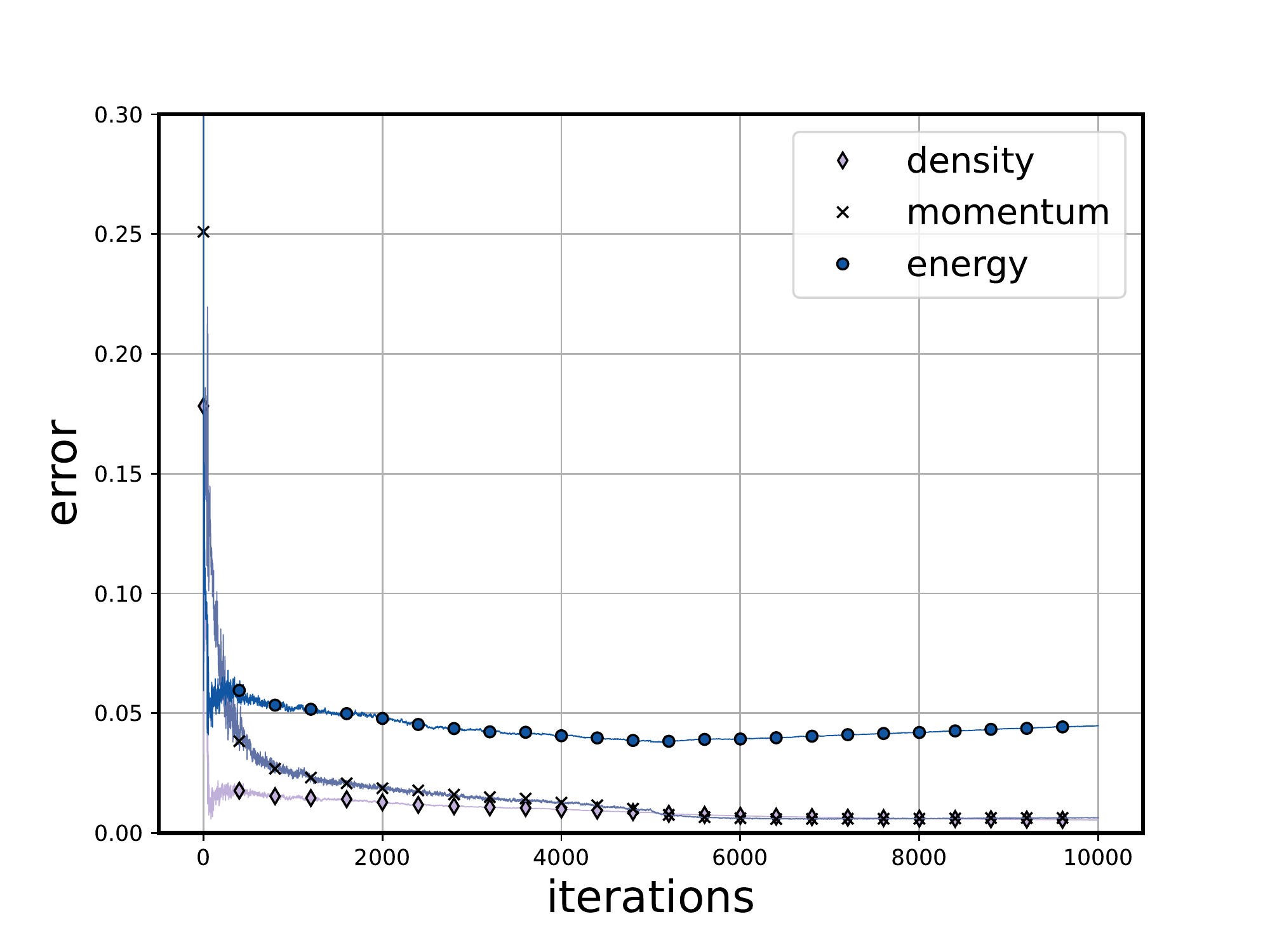}
    \includegraphics[width=0.45\textwidth]{ 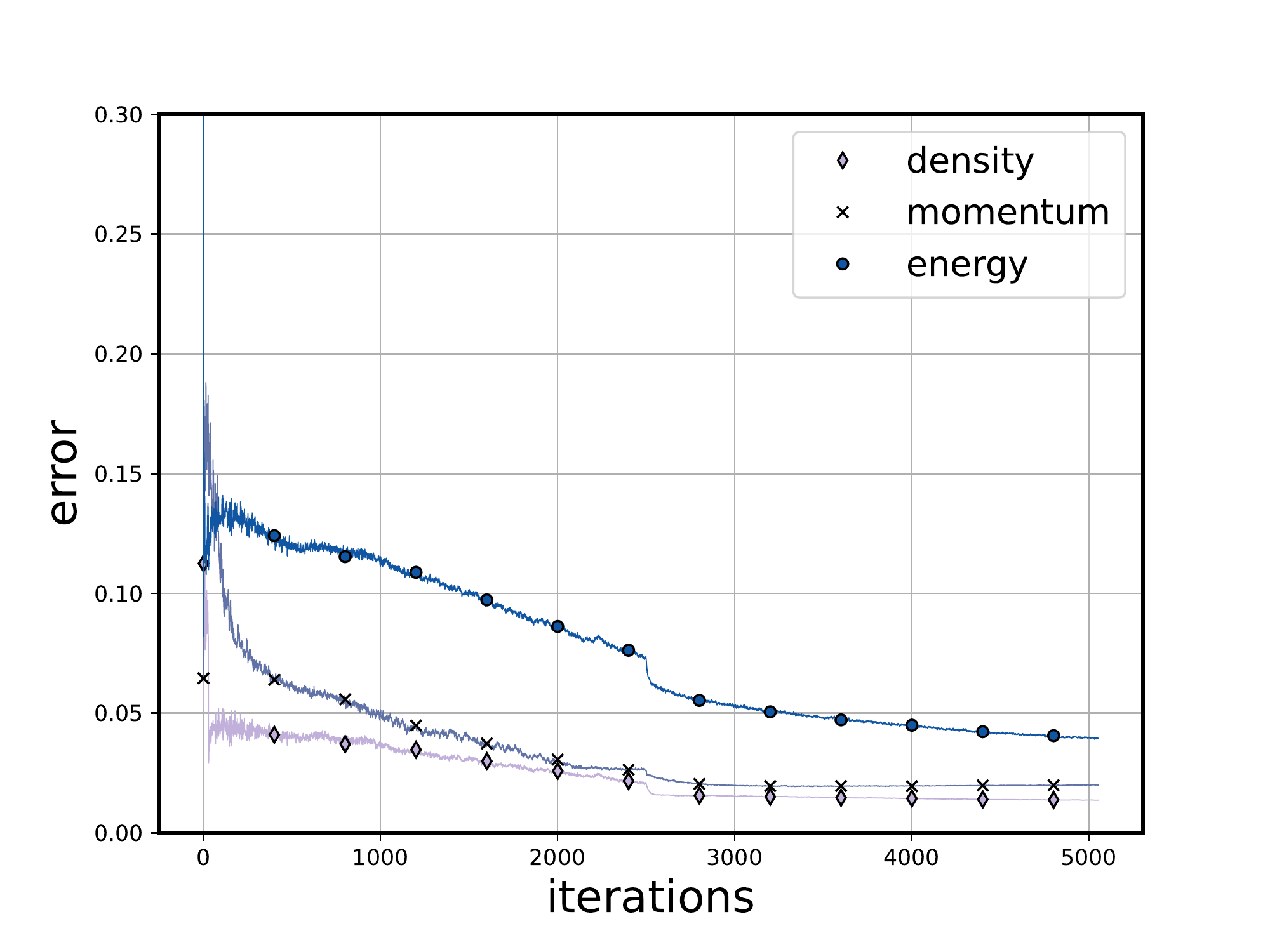}
    }

    \subfigure[{{\it Left:} Case III with $\varepsilon=10^{-3}$} and {{\it Right:} Case IV with  $\varepsilon=1$}.]
    {
    \includegraphics[width=0.45\textwidth]{ 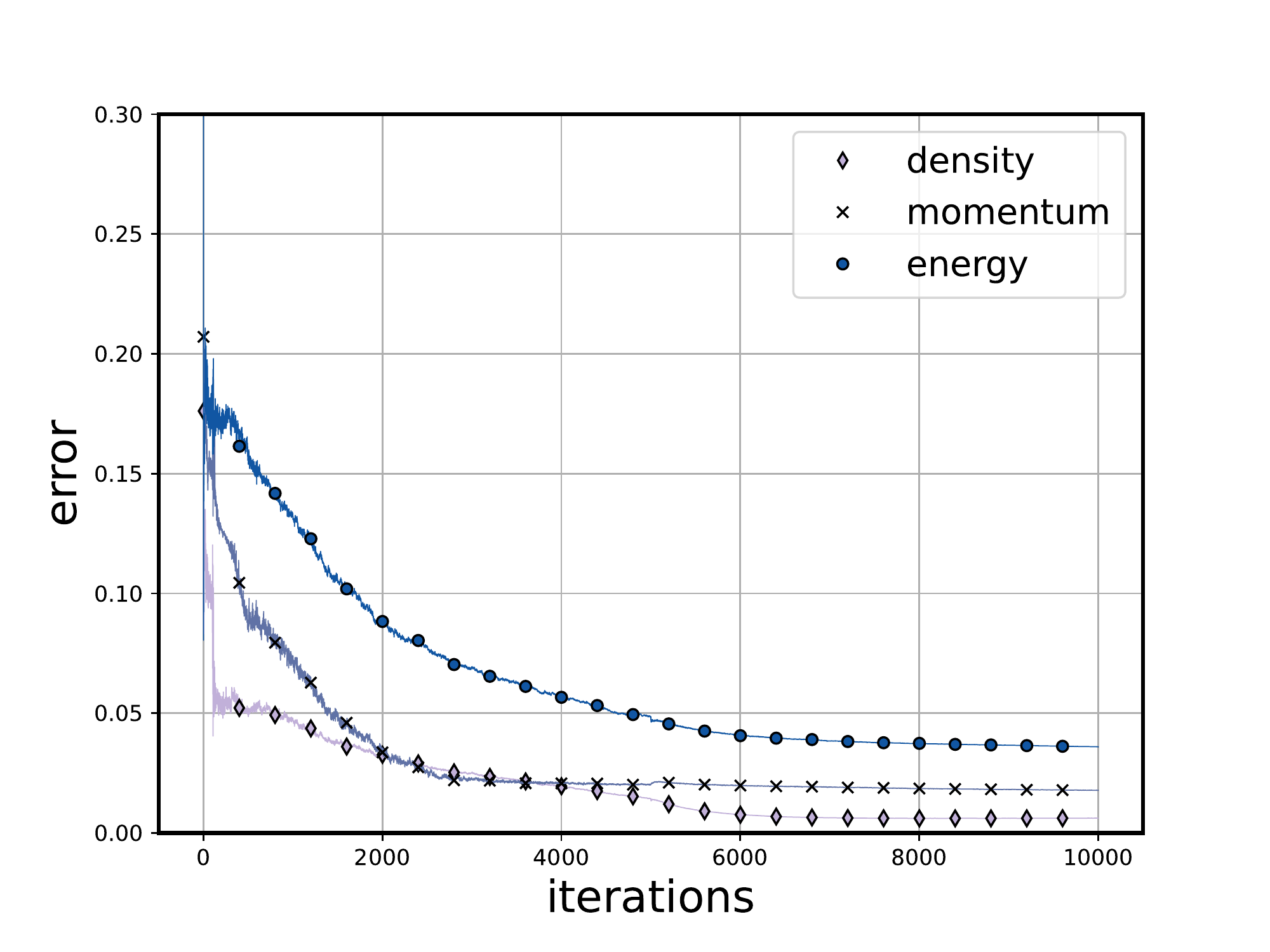}
    \includegraphics[width=0.45\textwidth]{ 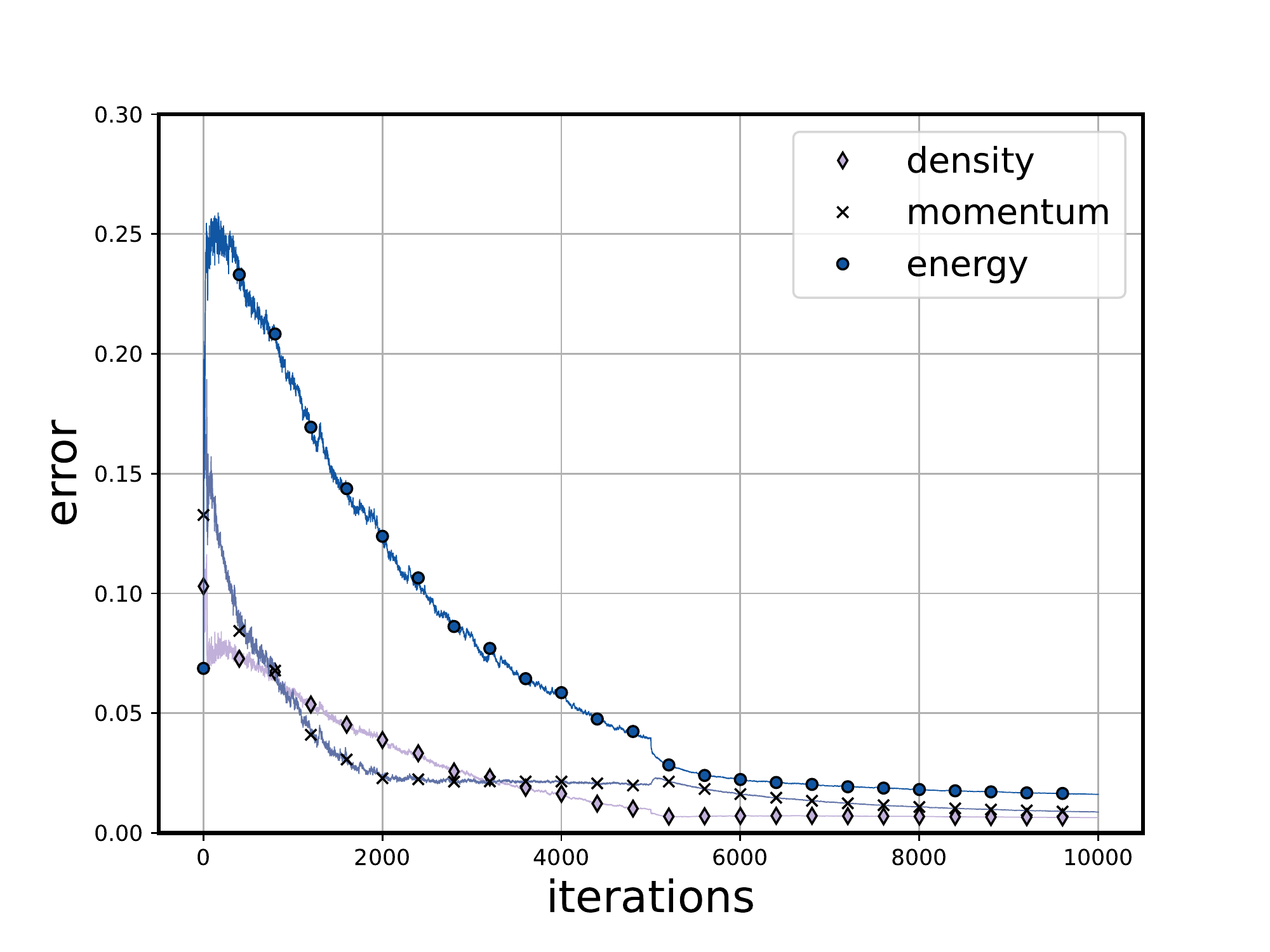}
    }
    \caption{Problem 2. Plot of the change of error for density, momentum and energy approximated by APNNs. }\label{fig:error-changes}
\end{figure}

Fig. \ref{fig:error-changes} illustrates the relative $l^2$ discrepancy in density, momentum, and energy throughout the training process.
An intriguing observation emerges, wherein the convergence of density and momentum outpaces that of energy across all instances, encompassing $\varepsilon = 1$ and $10^{-3}$.
Training the energy component proves notably more challenging.

Table~\ref{tab:losses} documents the mean square error of risks $\mathcal{R}^{\varepsilon}_{\text{residual}}, \mathcal{R}^{\varepsilon}_{\text{claw}}$, and $\mathcal{R}^{\varepsilon}_{\text{constraint}}$ across these scenarios.

\begin{small}
    \begin{table}[tbhp]
        \caption{Mean square errors of riskes $\mathcal{R}^{\varepsilon}_{\text{residual}}, \mathcal{R}^{\varepsilon}_{\text{claw}}$ and $\mathcal{R}^{\varepsilon}_{\text{constraint}}$ in terms of four cases.}\label{tab:losses}
        \centering
        \begin{tabular}{cccccccc}
            \toprule[1pt]
            \noalign{\smallskip}
            \multirow{2}*{\diagbox{{Case}}{Risk}}
             & \multicolumn{1}{c}{$\mathcal{R}^{\varepsilon}_{\text{residual}}$} & \multicolumn{3}{c}{$\mathcal{R}^{\varepsilon}_{\text{claw}}$} & \multicolumn{3}{c}{$\mathcal{R}^{\varepsilon}_{\text{constraint}}$}                                                                                                                        \\
             &                                                                   & \multicolumn{1}{c}{density}                                   & \multicolumn{1}{c}{momentum}                                        & \multicolumn{1}{c}{energy} & \multicolumn{1}{c}{0-order} & \multicolumn{1}{c}{1-order} & \multicolumn{1}{c}{2-order} \\
            \noalign{\smallskip}
            \midrule[1pt]
            \noalign{\smallskip}
            \multirow{1}*{{{I}}}
             & $3.92{\text{e-}6}$                                                & $7.01{\text{e-}7}$                                            & $1.70{\text{e-}5}$                                                  & $6.42{\text{e-}7}$         & $1.80{\text{e-}6}$          & $1.77{\text{e-}5}$          & $8.17{\text{e-}7}$          \\
            \multirow{1}*{{{II}}}
             & $1.30{\text{e-}5}$                                                & $1.16{\text{e-}5}$                                            & $1.97{\text{e-}5}$                                                  & $9.58{\text{e-}6}$         & $1.67{\text{e-}5}$          & $4.20{\text{e-}5}$          & $1.63{\text{e-}5}$          \\
            \multirow{1}*{{III}}
             & $1.54{\text{e-}5}$                                                & $1.81{\text{e-}6}$                                            & $4.67{\text{e-}6}$                                                  & $3.73{\text{e-}6}$         & $1.95{\text{e-}6}$          & $1.61{\text{e-}5}$          & $8.20{\text{e-}7}$          \\
            \multirow{1}*{{{IV}}}
             & $6.24{\text{e-}5}$                                                & $5.90{\text{e-}6}$                                            & $1.30{\text{e-}5}$                                                  & $2.09{\text{e-}5}$         & $3.62{\text{e-}6}$          & $6.08{\text{e-}6}$          & $5.07{\text{e-}6}$          \\
            \noalign{\smallskip}
            \bottomrule[1pt]
        \end{tabular}
    \end{table}
\end{small}

{\bf{Case V:}} $\varepsilon=10^{-2}$
\begin{equation*}
    \begin{aligned}
         & \rho_0(x) =
        \left \{
        \begin{aligned}
             & 1.5, \; x < 0,    \\
             & 0.625 \; x \ge 0,
        \end{aligned}
        \right.,       \\
         & u_0(x) =0,  \\
         & T_0(x) =
        \left \{
        \begin{aligned}
             & 1.5, \; x < 0,   \\
             & 0.75 \; x \ge 0.
        \end{aligned}
        \right.
        .
    \end{aligned}
\end{equation*}

\begin{figure}[ht]
    {
        \includegraphics[width=0.45\textwidth]{ 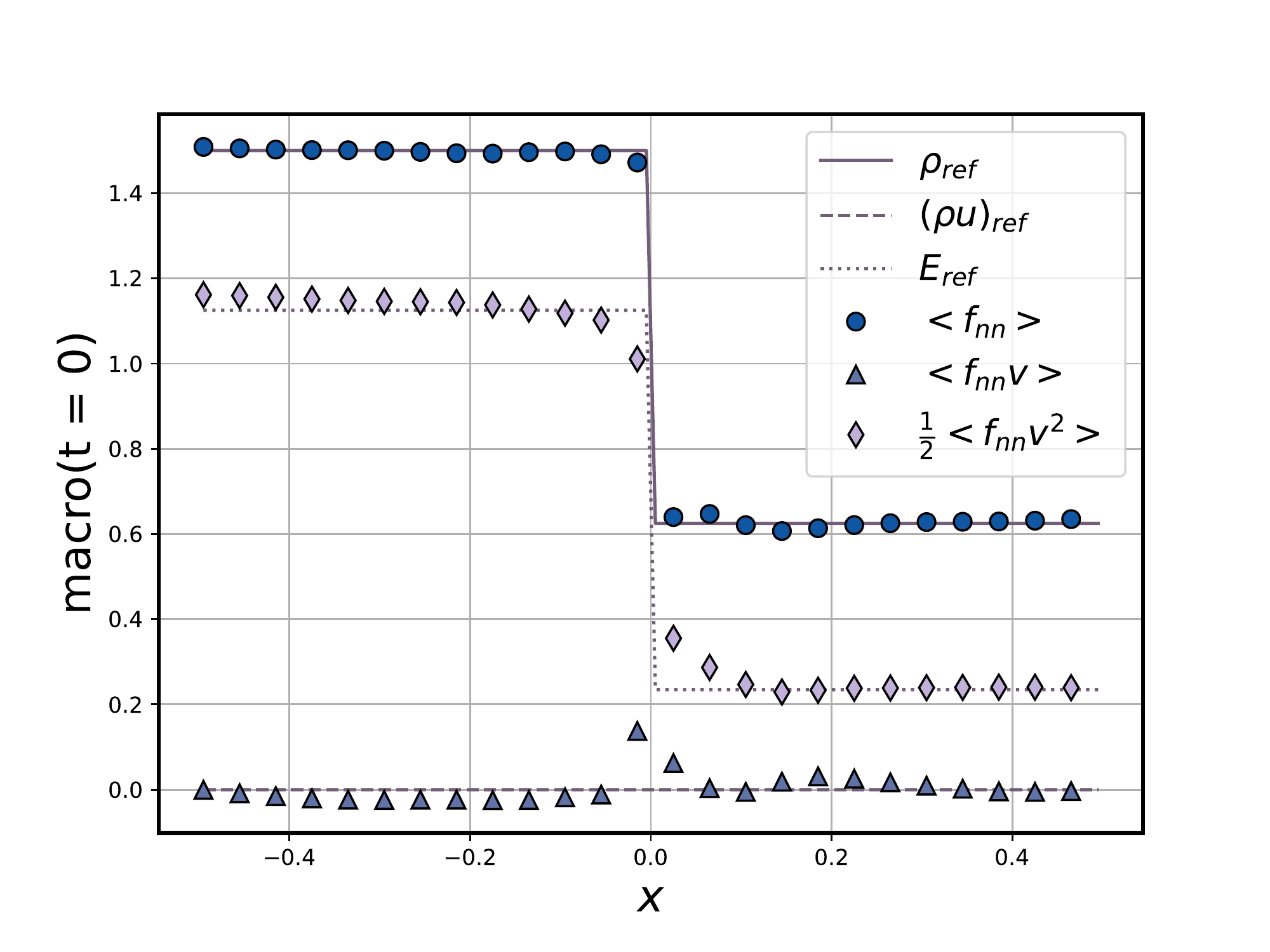}
        \includegraphics[width=0.45\textwidth]{ 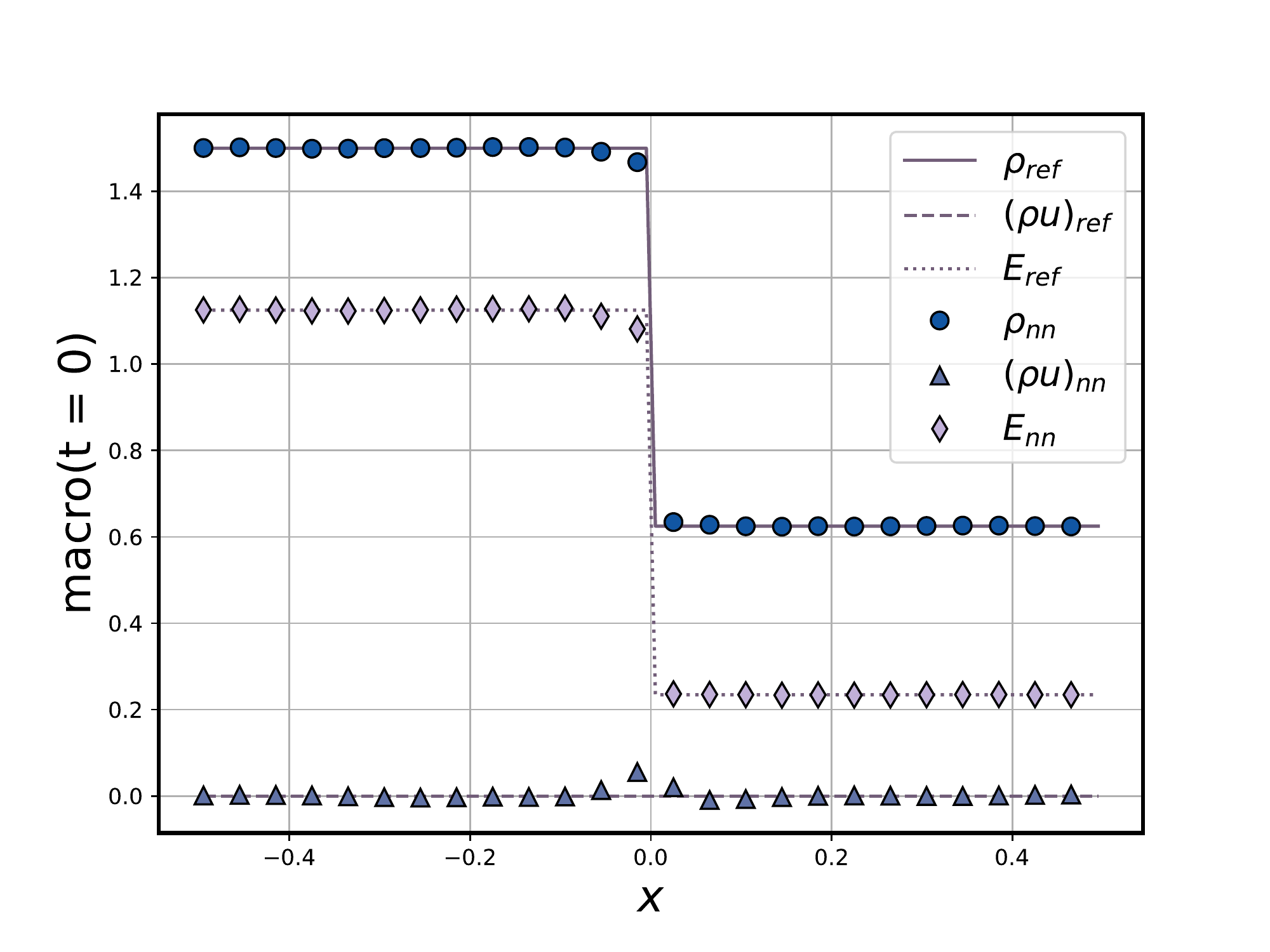}
    }
    \caption{Problem 2---Case V. Plot of density, momentum and energy at time $t = 0$: Approximated by APNNs (marker) vs. Ref (line).
        $\varepsilon = 10^{-2}$ and the units of neural networks are $[3, 128, 128, 128, 128, 128, 128, 1]$ for $f$ and $[2, 96, 96, 96, 96, 96, 96, 1]$ both for $\rho, u$ and $T$.
        Batch size is $512$ in domain, and $512$ on initial condition. ${\lambda_{3}} = (10, 10, 1), {{\lambda_5}} = (1, 100, 100), \lambda_6 = 100$ and others are set to be $1$.}\label{fig:example-discontinuous}
\end{figure}
The figure depicted in Fig. \ref{fig:example-discontinuous} illustrates the graphical representation of the estimated macroscopic variables at the instant of $t = 0$.

\begin{figure}[ht]
    \subfigure[{The approximate $\rho_\text{nn}, u_\text{nn}, T_\text{nn}$} vs.\ reference solutions. {\it Left: $t = 0.005$} and {\it Right: $t = 0.01$}.]
    {
        \includegraphics[width=0.45\textwidth]{ 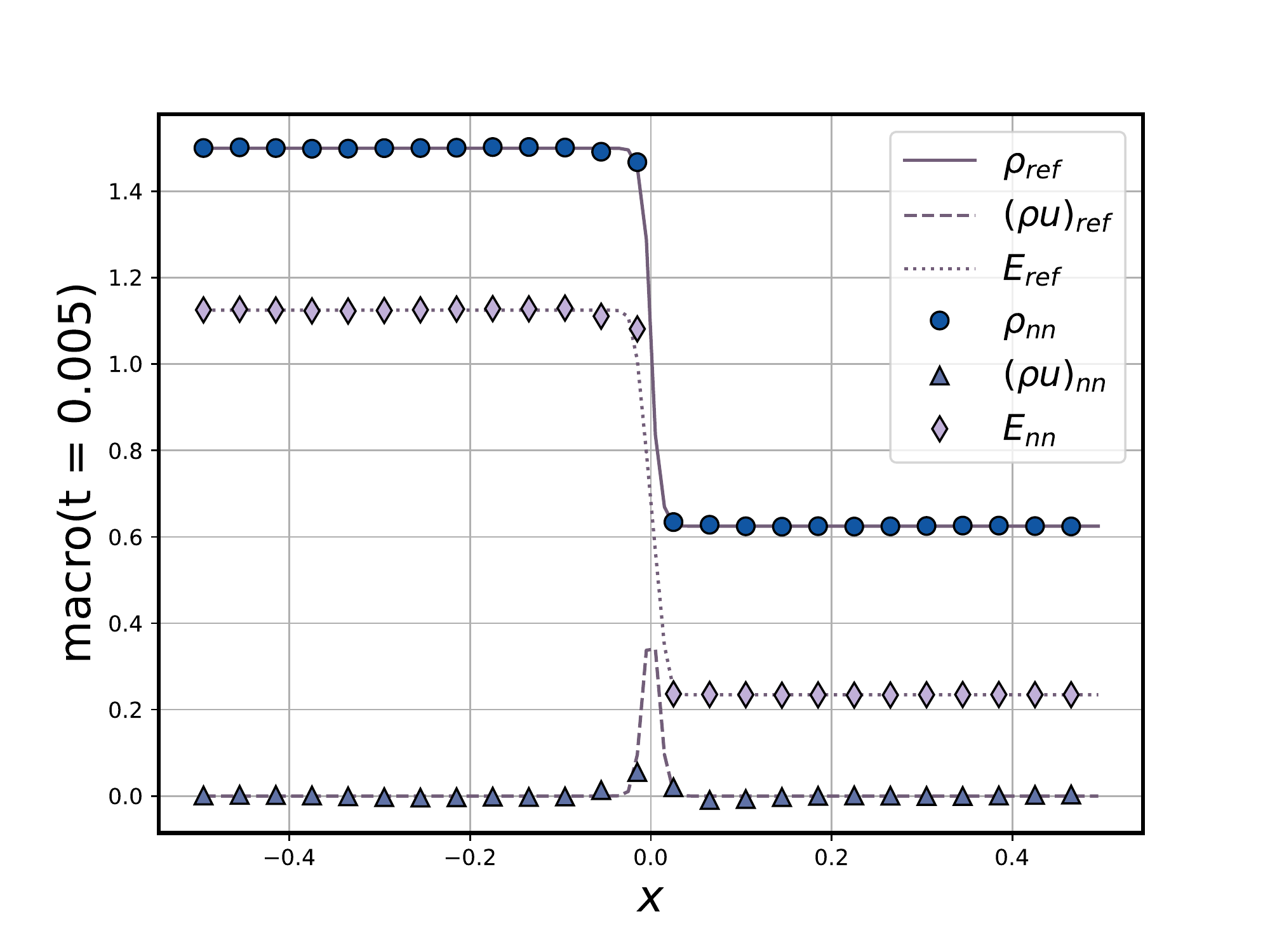}
        \includegraphics[width=0.45\textwidth]{ 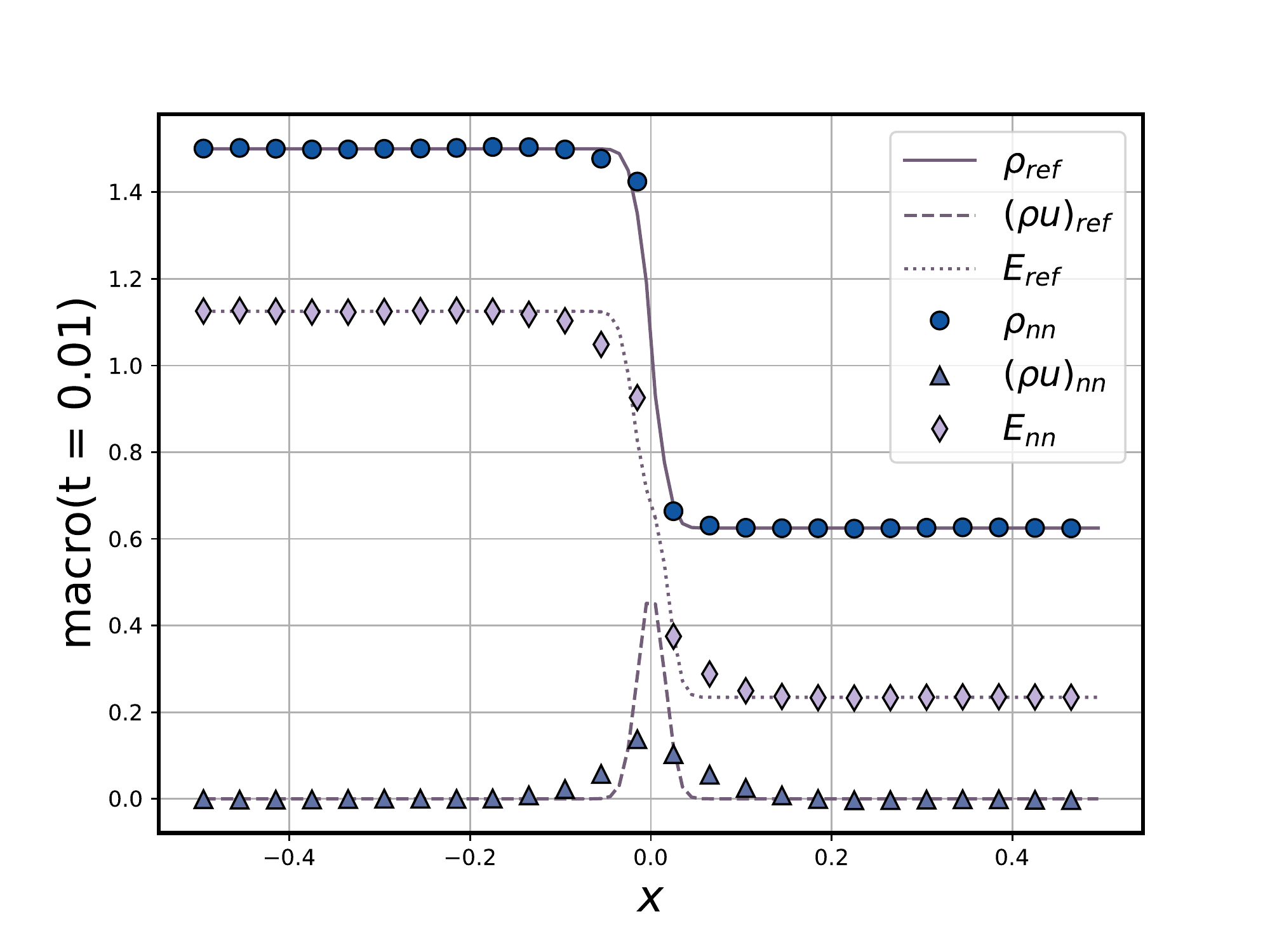}
    }
    \subfigure[{The approximate $\rho_\text{nn}, u_\text{nn}, T_\text{nn}$} vs.\ reference solutions. {\it Left: $t = 0.02$} and {\it Right: $t = 0.03$}.]
    {
        \includegraphics[width=0.45\textwidth]{ 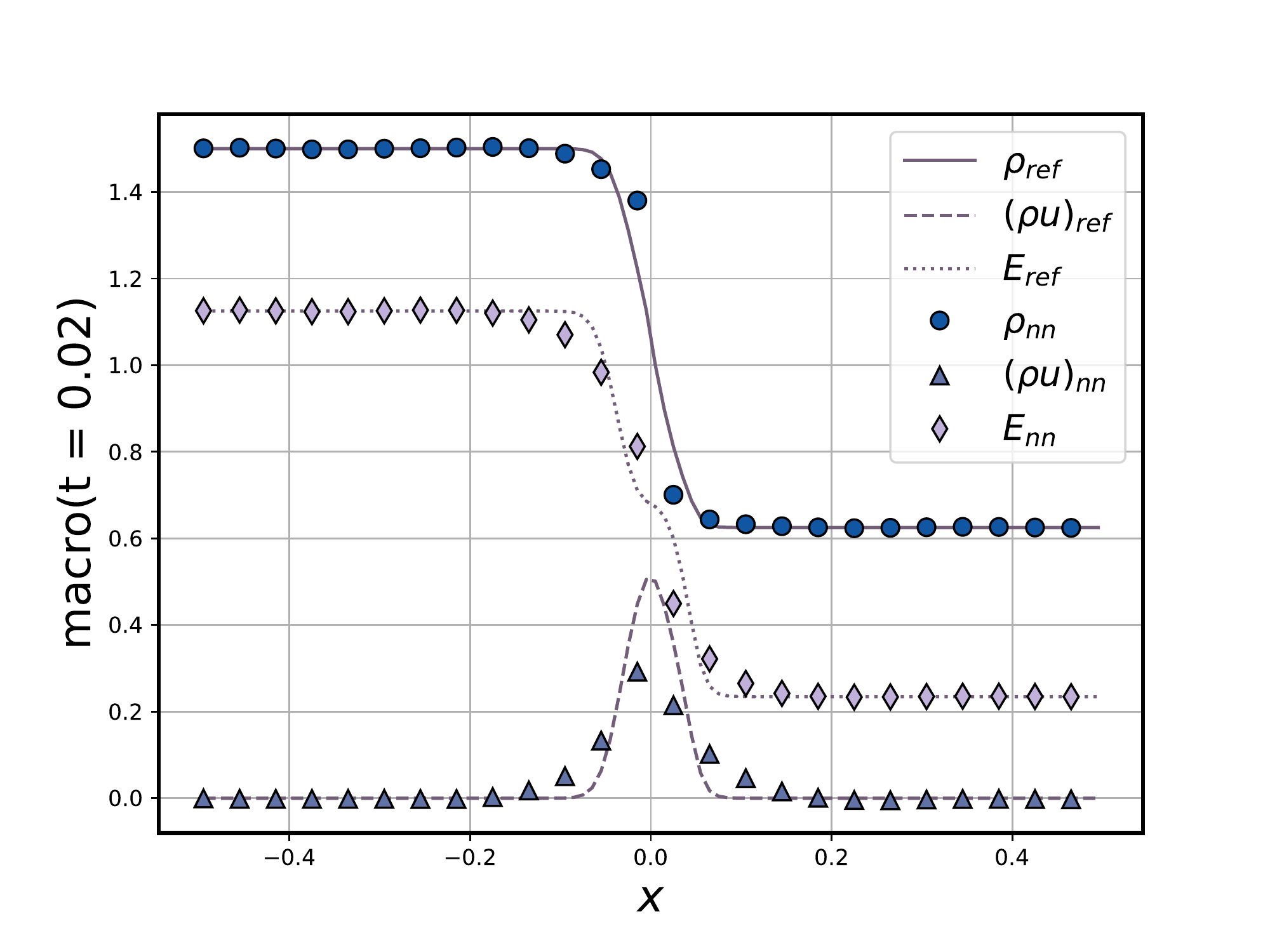}
        \includegraphics[width=0.45\textwidth]{ 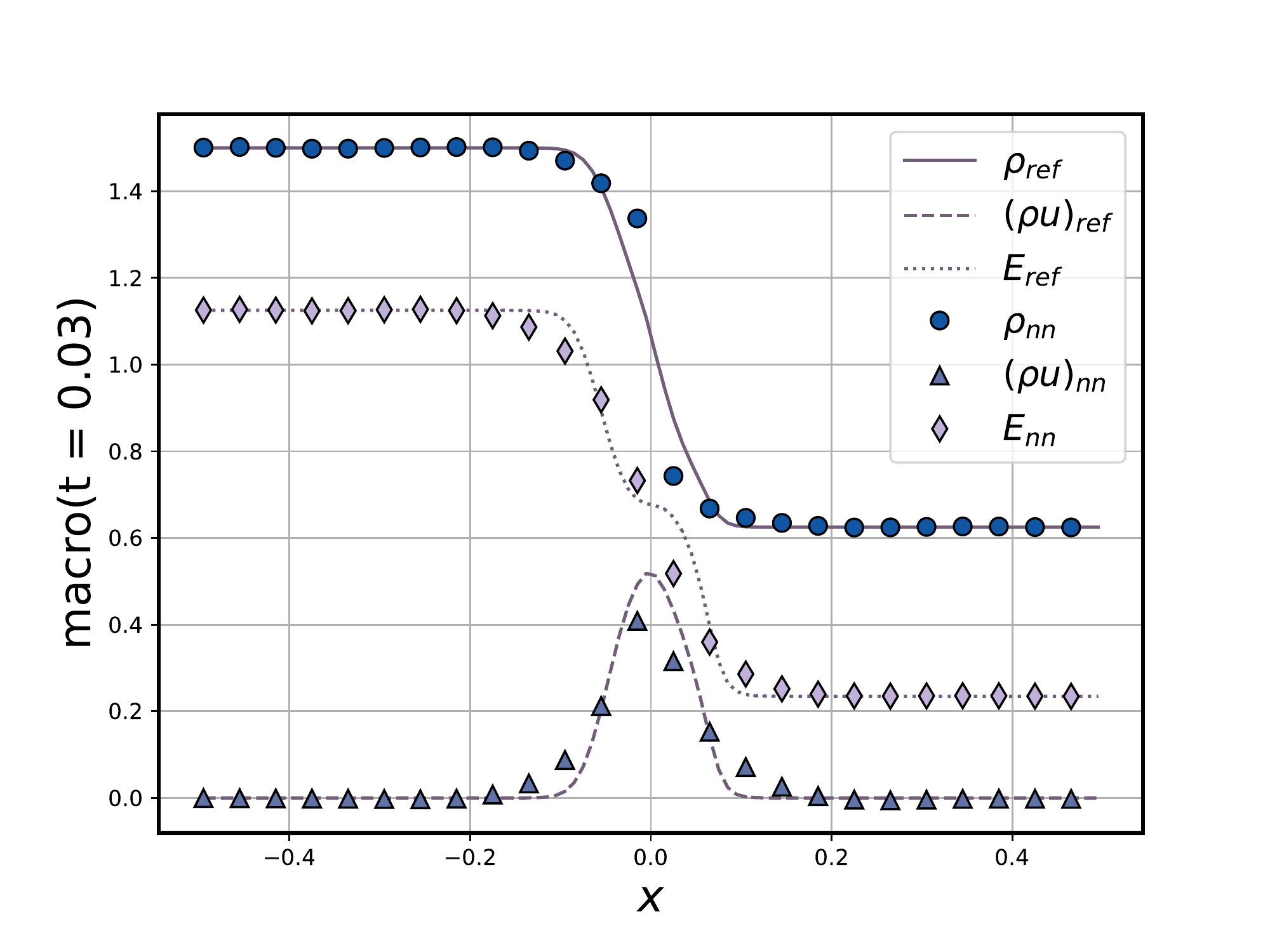}
    }
    \subfigure[{The approximate $\rho_\text{nn}, u_\text{nn}, T_\text{nn}$} vs.\ reference solutions. {\it Left: $t = 0.04$} and {\it Right: $t = 0.05$}.]
    {
        \includegraphics[width=0.45\textwidth]{ 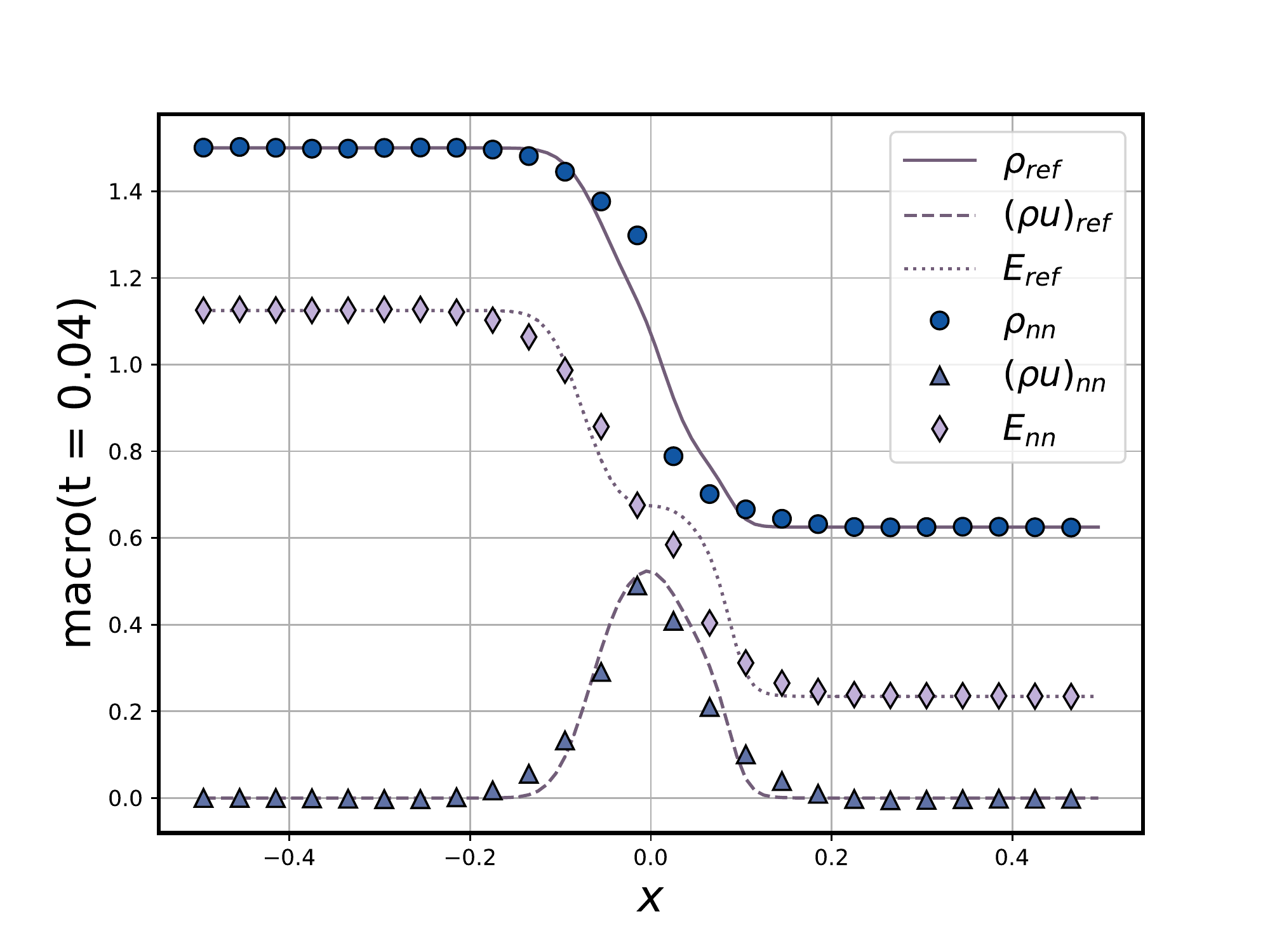}
        \includegraphics[width=0.45\textwidth]{ 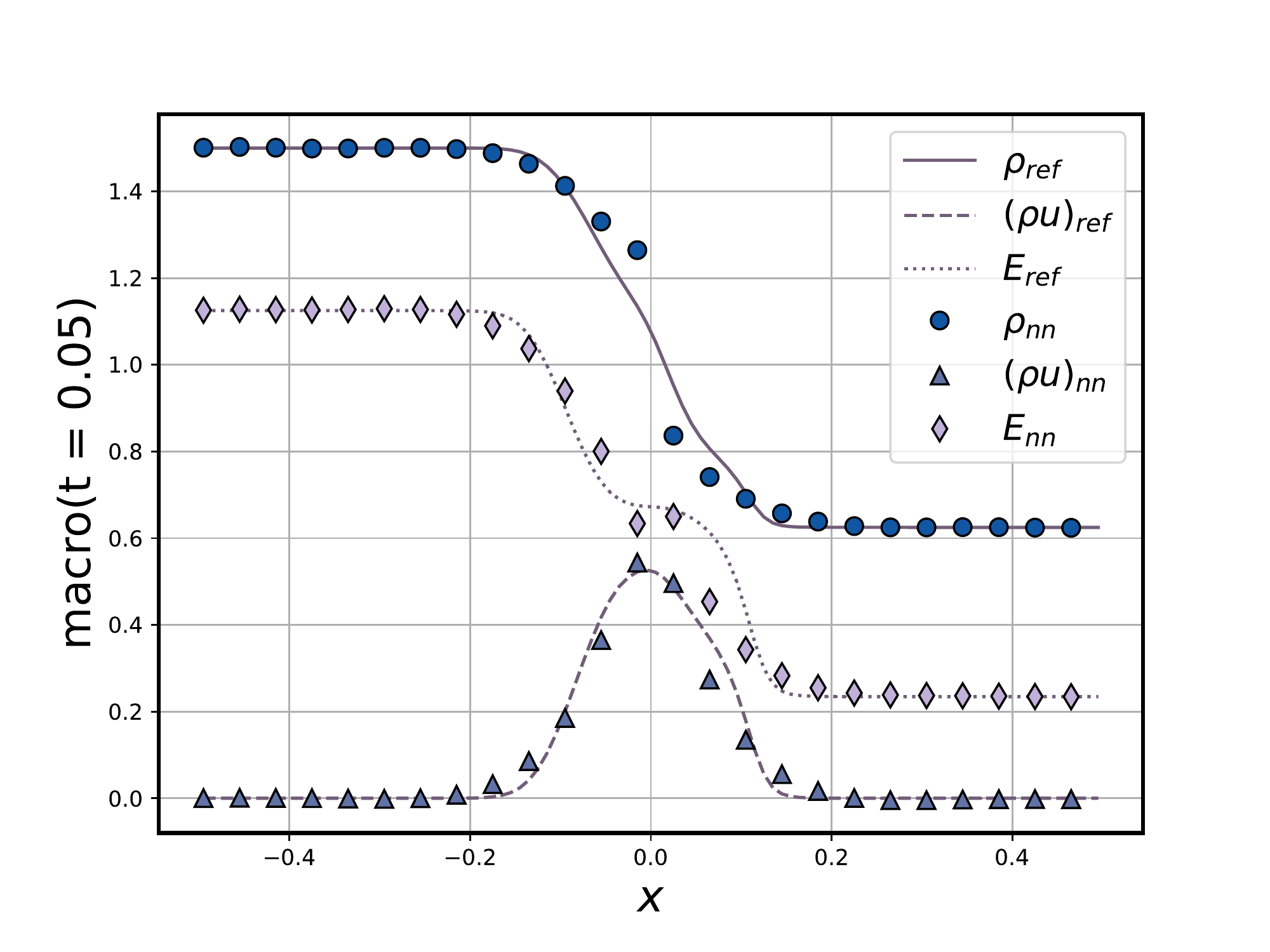}
    }
    \caption{Plot of density, momentum and energy at time $t = 0.005$ to $0.05$: Approximated by APNNs (marker) vs. Ref (line).}\label{fig:example-disprocess}
\end{figure}
Fig. \ref{fig:example-disprocess} depicts the graphical representation of the approximated macroscopic quantities spanning from $t = 0.005$ to $t = 0.05$.
Notably, at $t = 0.005 < \varepsilon$, a discernible transition in the macroscopic solutions becomes evident.

\section{Discussion and conclusion}

In this paper, we have devised novel and efficient Asymptotic-Preserving Neural Networks (APNNs) to numerically approximate the multiscale kinetic equations involving possibly small Knudsen number.
For the linear transport equation with diffusive scaling, our APNN, utilize parity equations, by designing the loss function that captures physical conservation property, and also the inflow boundary condition.
For the BGK equation, with the a collision term that exhibits both nonlinearity and non-locality, along with a boundary layer effect,  we have devised an APNN method that with a loss that not only follows the equation but also the dynamics of the conserved moments,  in addition to precisely enforces the boundary conditions.

During the training process, a phenomenon becomes apparent: the convergence of momentum and energy is slower in comparison to that of density. This discrepancy is likely influenced by the value of $\varepsilon$. The question that arises is how to explain this phenomenon, which needs further study.

Due to the nonlinearity, high dimensionality and multiscale nature of kinetic equations, the construction of the neural network holds tremendous tremendous. Our studies are preliminary but call for more investigations to more physical equations such as the Boltzmann and Landau equations.

\section*{Acknowledgement}
This work is partially supported by the National Key R\&D Program of China Project No. 2020YFA0712000
and Shanghai Municipal of Science and Technology Major Project No. 2021SHZDZX0102. Shi Jin
is also supported by NSFC grant No. 11871297. Zheng Ma is also supported by NSFC Grant
No. 12031013, No.92270120 and partially supported by Institute of Modern Analysis---A Shanghai Frontier Research Center.


\end{document}